\documentclass[11pt,oneside,reqno]{amsart}
\usepackage{amsfonts}
\usepackage{amsmath}
\usepackage{color}
\usepackage{amssymb}
\usepackage{multirow}
\usepackage{lscape}
\input xy \xyoption{all}
\newcommand\fg{{\frak{g}}} 
\newcommand\fh{{\frak h}}

\newcommand\fm{{\frak m}} 
\newcommand\fn{{\frak n}} 
\newcommand\fp{{\frak p}} 
\newcommand\fa{{\frak a}} 
\newcommand\fb{{\frak b}} 
\newcommand\fd{{\frak d}}

\newcommand\fs{{\frak s}} 
 
\newcommand\fu{{\frak u}} 
\newcommand\fv{{\frak v}}
 
\newcommand\fy{{\frak y}}

\newcommand{\fgl}{\mathop{{\frak g \frak l}}} 
\newcommand{\fsl}{\mathop{{\frak s \frak l}}}

\newcommand{\fso}{\mathop{{\frak s \frak o}}}

\newcommand\RR{{\mathbb R}}

\newcommand\PP{{\mathbb P}}
\newcommand\ph{\varphi} 

\newtheorem{theo}{Theorem}[section] 
\newtheorem{pr}[theo]{Proposition} 
\newtheorem{de}[theo]{Definition} 
\newtheorem{ex}[theo]{Example} 
 
\newtheorem{co}[theo]{Corollary} 
\newtheorem{lm}[theo]{Lemma} 
\newcommand{\GL}{\mathop{{\rm GL}}} 
\newcommand{\SO}{\mathop{{\rm SO}}}

\newcommand{\SL}{\mathop{{\rm SL}}}

\newcommand{\Spin}{\mathop{{\rm Spin}}} 
\newcommand{\G}{\mathop{{{\rm G}_{2}^*}}} 

\newcommand{\End}{\mathop{{\rm End}}}
\newcommand{\Ad}{\mathop{{\rm Ad}}} 
\newcommand\ip{{\langle\cdot \,,\cdot \rangle}} 
\newcommand\la{\langle}
\newcommand\ra{\rangle}
\newcommand{\Span}{\mathop{{\rm span}}} 
\newcommand{\tr}{\mathop{{\rm tr}}} 
\newcommand{\diag}{\mathop{{\rm diag}}} 
\newcommand{\sgn}{\mathop{{\rm sgn}}} 
\begin{document}

\title[Holonomy groups of $G_2^*$-manifolds]{Holonomy groups of $G_2^*$-manifolds}
\date{\today}
\author{Anna Fino and Ines Kath}
\maketitle
\begin{abstract}
We classify the holonomy algebras of manifolds admitting an indecomposable torsion free $G_2^*$-structure, i.e. for which the holonomy representation does not leave invariant any proper non-degenerate subspace. We realize some of these Lie algebras as holonomy algebras of left-invariant metrics on Lie groups.
 \end{abstract}
\section{Introduction}
Holonomy groups are a useful tool in the study of semi-Riemannian manifolds. They make it possible to apply algebraic methods to geometric problems such as the existence of special geometric structures or the decomposability of manifolds. So it is natural to ask which Lie groups can be the holonomy group of a semi-Riemannian manifold. Since we are mainly interested in connected holonomy groups we may equivalently ask which Lie subalgebras of $\fso(p,q)$ are holonomy algebras. For Riemannian manifolds, there is a complete answer to this question. Berger's list gives a classification of irreducible Riemannian holonomy algebras of non-locally symmetric spaces \cite{Be, Br}. Moreover, holonomy algebras of locally symmetric Riemannian manifolds can be read off from Cartan's classification of Riemannian symmetric spaces. The pseudo-Riemannian situation is much more complicated. In general, a holonomy representation, i.e., the natural representation of a holonomy group on the tangent space can have isotropic invariant subspaces and is not necessarily completely reducible. Therefore it does not suffice to determine all irreducible holonomy groups. A complete classification is only known for Lorentzian manifolds, it is due to Leistner \cite{Le}. For metrics of index greater than one only partial results are known. For instance, there are results for manifolds with special geometric structure. Galaev classified holonomy algebras of pseudo-K\"ahlerian manifolds of index~2 \cite{Ga}. Furthermore, holonomy groups of pseudo-quaternionic-K\"ahlerian manifolds of non-zero scalar curvature were classified by Bezvitnaya \cite{Bz}.

In the present paper, we want to turn to another special geometry, namely the pseudo-Riemannian analogue of a torsion-free $G_2$-structure, which is well known from the holonomy theory of Riemannian manifolds since $G_2$ is one of the groups on Berger's list. While torsion-free $G_2$-structures exist on Riemannian 7-manifolds, their pseudo-Riemannian analogues are structures on manifolds of signature $(4,3)$. They are characterised by the fact that their holonomy is contained in the non-compact subgroup $G^*_2\subset SO(4,3)$  of type $G_2$, which is defined as the stabiliser of a certain generic 3-form.
There are other nice characterisations of this group, e.g., $\G$ is the stabiliser of a non-isotropic element of the real spinor representation of $\Spin(4,3)$ and it can also be understood as the stabiliser of a cross product on $\RR^{4,3}$. Hence a  torsion-free $G^*_2$-structure on a pseudo-Riemannian manifold $M$ of signature $(4,3)$ can be understood as a parallel generic 3-form, a parallel non-isotropic spinor field or a parallel cross-product `$\times$' on $M$.

Our aim is to classify the holonomy algebras of manifolds admitting a torsion-free  $\G$-structure, where we want to assume that this $\G$-structure is indecomposable, that is, its holonomy representation does not leave invariant any proper non-degenerate subspace. By a classification we mean a classification as subalgebras of $\fg^*_2\subset\fso(4,3)$ up to conjugation by elements of $\SO(4,3)$.

There are already some results in this direction. In \cite{Ka} indefinite symmetric spaces with $\G$-structure are classified. Their holonomy  algebras can be read off from this classification. It turns out that they are abelian and of dimension two or three. 
Furthermore, some results on left-invariant torsion-free $\G$-structures on Lie groups (or, equivalently, $\G$-structures on Lie algebras) are known. Examples of (decomposable) torsion-free $\G$-structures with $1$-dimensional and
$2$-dimensional holonomy have been found by M. Freibert \cite{Fr}  on  almost abelian
Lie algebras. In \cite{FL}   Fino and Lujan studied torsion-free  $\G$-structures with holonomy  $\G$ on nilpotent Lie algebras, showing in particular that, up to isomorphism, there exists only one indecomposable nilpotent Lie algebra admitting a torsion-free  $\G$--structure such that the center is definite with respect to the induced inner product.  In \cite{FL} an example of an indecomposable  torsion-free  $\G$-structure with $6$-dimensional holonomy   on a nilpotent Lie algebra is also given.  

The first step in the classification of holonomy algebras is to get algebraic conditions for candidates for holonomy algebras $\fh$ strictly contained in $\fg_2^*$.  These conditions can be derived from the following three facts. Firstly, since $\fh$ is a proper subalgebra of $\fg_2^*$, the natural representation of $\fh\subset \fg_2^*\subset \fso(4,3)$ on $\RR^{4,3}$ has to leave invariant an isotropic subspace. Secondly, several restrictions come from the indecomposability of this representation. Thirdly, every holonomy algebra is a Berger algebra, i.e., it satisfies Berger's first criterion, which gives further conditions for $\fh$.  In this paper we give a complete answer to this algebraic part of the classification problem. That is, we classify indecomposable Berger algebras strictly contained in $\fg_2^*$.  The results can be summarised as follows. Let $\fp_1$, $\fp_2$ denote the two 9-dimensional parabolic subalgebras of $\fg_2^*$, which can be characterised by the action of $G_2^*$ on isotropic subspaces of $\RR^{4,3}$. The action of $G_2^*$ on isotropic lines is transitive and $\fp_1$ is the Lie algebra of the stabiliser of an isotropic line.  Furthermore, the action of $G_2^*$ on 2-planes $E=\Span\{b_1,b_2\}$ satisfying $b_1\times b_2=0$ is transitive and $\fp_2$ is the Lie algebra of the stabiliser of such a 2-plane. We have $\fp_1\cong\fgl(2,\RR)\ltimes \fm$, where $\fm$ is three-step nilpotent and $\fp_2\cong\fgl(2,\RR)\ltimes \fn$, where $\fn$ is two-step nilpotent. Both $\fp_1$ and $\fp_2$ are indecomposable Berger algebras. We will distinguish arbitrary indecomposable Berger algebras $\fh\subset\fg_2^*\subset \fso(4,3)$ by the dimension of the socle of their natural representation on $\RR^{4,3}$. The socle is the maximal semisimple subrepresentation. By indecomposability, it is isotropic. We will  say that $\fh$ is of Type I, II or III, if the dimension of the socle is one, two or three. In particular, $\fp_1$ is of Type I, $\fp_2$ is of Type II. We show that $\fh\subset \fp_1$ up to conjugation if $\fh$ is of Type I or III. If $\fh$ is of Type II, then $\fh\subset \fp_2$ up to conjugation. Let $\fa$ be the projection of $\fh$ to $\fgl(2,\RR)\subset\fp_i$ for $i=1,2$, respectively. Then we may assume that $\fa$ is one of the representatives of conjugacy classes of subalgebras of $\fgl(2,\RR)$. Roughly speaking, for each of these representatives we classify the subalgebras $\hat \fm\subset \fm$ and $\hat\fn\subset\fn$ for which $\fa\ltimes \hat\fm$ and $\fa\ltimes \fn$ are indecomposable Berger algebras. Finally, for each type, we get a list of all (conjugacy classes of) indecomposable Berger algebras
(Theorems \ref{T1}, \ref{T3} and \ref{T2}).

The second part of the classification consists in the realisation of the possible holonomy algebras by metrics. As we already mentioned above,  left-invariant metrics on 7-dimensional nilpotent and solvable Lie groups may provide interesting examples  of such metrics. In Section~\ref{S4}, we  give new examples of left-invariant metrics with holonomy contained in $\fg_2^*$. In particular, we can  provide examples for each of the Types I, II and III. As for Type I, we can realise $\fm$ and, furthermore, a  7-dimensional solvable Lie algebra and a 6-dimensional nilpotent one as a holonomy algebra. Besides $\fn$ and $\fsl(2,\RR)\ltimes \fn$, we give a 3-dimensional abelian example of Type II. Finally, we can realise a three-dimensional abelian Lie algebra of Type III. 
Another special class of pseudo-Riemannian manifolds is that of symmetric spaces. As already mentioned, symmetric spaces with $\G$-structure were determined in \cite{Ka}. In Section~\ref{S3}, we check how their holonomy algebras fit into the classification. 

\textsc{Acknoledgement} We are very thankful for the opportunity to stay two weeks in the very stimulating atmosphere at Mathematisches Forschungsinstitut Oberwolfach, where we finished essential parts of this paper. 

\textsc{Notation} If $b_1,\dots,b_n$ is a basis of a vector space $W$, then we denote by $b^1,\dots,b^n$ its dual basis of $W^*$. Furthermore, $b^{i_1\dots\, i_k}:=b^{i_1}\wedge\dots\wedge b^{i_k}\in\bigwedge^k W^*$, $b_{i_1\dots\, i_k}:=b_{i_1}\wedge\dots\wedge b_{i_k}\in\bigwedge^k W$ and $b^j_i:=b_i\otimes b^j\in W\otimes W^*\cong \End(W)$.
\section{ Holonomy groups contained in $G^*_2$ } \label{S2}
\subsection{The group $\G$}
Let $M$ be a simply connected manifold of signature $(4,3)$. Suppose that $M$ admits a $\G$-structure, given by a generic three-form $\omega$. Then $\omega$ defines an orientation and a metric $g$ of signature $(4,3)$ on $M$. Here signature $(4,3)$ means that $g=\diag(-1,-1,-1,-1,1,1,1)$ with respect to a suitable basis. Let $\omega$ be parallel with respect to $g$.  Then the holonomy group $H$ of $(M,g)$ is contained in  $G_2^*$. 

Equivalently, we could have started with a pseudo-Riemannian manifold $(M,g)$ of signature $(4,3)$ whose holonomy group $H$ is contained in $\G\subset \SO(4,3)$. Then $\G$ defines a parallel 3-form $\omega$ on $(M,g)$ and $g$ is induced by~$\omega$.

Let $\fh$ denote the Lie algebra of $H$. Suppose that the holonomy representation of $\fh$ on $V:=T_x(M)$ is indecomposable but not irreducible.

 A subspace $0\not=U\subset V$ is called isotropic if $g(u,u)=0$ for all $u\in U$.
 
Let us give explicit formulas. We choose a basis $e_1,\dots,e_7$ of $V$ such that the 3-form $\omega$ equals
$$
\omega_0=\sqrt 2(e^{167}+e^{235})-e^4\wedge (e^{15}-e^{26}-e^{37})\,.
$$

Then the orientation of $V$ is defined to be the orientation of $e_1,\dots,e_7$ and the induced metric equals $\ip=2(e^1\cdot e^5+e^2\cdot e^6+e^3\cdot e^7)- (e^4)^2$. In particular, we can identify $\G$ with the subgroup of $\GL(7, {\mathbb R})$ that stabilises $\omega_0$. Then $\G\subset SO(4,3)$ with respect to $\ip$. The Lie algebra $\fg_2^*$ consists of all matrices of the form
\begin{equation}\label{EG2}
\left(
\begin{array}{ccccccc}
s_1+s_4 &-s_{10}&s_9&\sqrt 2 s_6&0&-s_{11}&-s_{12}\\
-s_8&s_1&s_2&\sqrt 2 s_9&s_{11}&0&s_6\\
s_7&s_3&s_4&\sqrt 2 s_{10}&s_{12}&-s_6&0\\
\sqrt2 s_5&\sqrt 2 s_7&\sqrt 2 s_8&0&\sqrt 2 s_6&\sqrt 2 s_9&\sqrt 2 s_{10}\\
0&s_{13}&s_{14}&\sqrt 2 s_5&-s_1-s_4&s_8&-s_7\\
-s_{13}&0&-s_5&\sqrt 2 s_7&s_{10}&-s_1&-s_3\\
-s_{14}&s_5&0&\sqrt 2 s_8&-s_9&-s_2&-s_4
	\end{array}\right),
\end{equation}	
where $s_1,\dots,s_{14}\in\RR$.

The 3-form $\omega$ defines a cross product on $V$ by
$$\la u\times v ,w\ra=\omega(u,v,w).$$
This cross product is antisymmetric and satisfies 
$$\la u\times v,u\ra=0,\quad u\times (u\times v)=-\la u,u\ra v+\la u,v\ra u,$$
for all $u,v\in V$.

The group $\G$ can also be understood as the stabiliser of a non-isotropic spinor $\psi_0$ in the real spinor representation $\Delta$ of $\Spin(V)\cong \Spin(4,3)$.  Indeed, the two-fold covering map $\lambda:\Spin(V)\rightarrow \SO(V)$ induces an isomorphism from this stabiliser to $\G$. This is well known, see e.g. \cite{Ka0} for details. We give also \cite{Ka} as a reference here since the above mentioned formulas for $\omega_0$ and $\fg_2^*$ can be obtained from the description of the Clifford algebra of $V$ in \cite{Ka} in replacing the basis $e_1,\dots,e_7$ used there by $e_7, e_5, e_6, e_4, e_3, e_1, e_2$. We want to recall the following well-known facts. The spinor module $\Delta$ admits an inner product $\ip_\Delta$ of signature $(4,4)$ that is invariant with respect to the Clifford multiplication, i.e.,  
$$\langle v \cdot \ph,\psi\rangle_\Delta+\langle \ph,v\cdot \psi\rangle_\Delta=0,$$
for all $v\in V$. The Clifford multiplication of the non-isotropic spinor $\psi_0$ by a vector
$$
V\longrightarrow \{\psi_0\}^\perp\subset \Delta,\quad v\longmapsto v\cdot\psi_0
$$
is an isomorphism from $V$ to $\{\psi_0\}^\perp$. The spinor $\psi_0$ defining $\G$ is related to the cross product by
\begin{equation}\label{Ecross}
u\cdot v\cdot\psi_0+\la u,v\ra\psi_0=(u\times v)\cdot \psi_0.
\end{equation}
The map 
$$
\Delta\ni \ph\longmapsto U(\ph):= \{ v\in V\mid v\cdot \ph=0\}\subset V
$$
induces a bijection from the set  of projective isotropic spinors to the set of  3-dimensional isotropic subspaces of $V$. 

\subsection{The type of a holonomy algebra contained in $\fg_2^*$}

Since the holonomy representation of $\fh$ is indecomposable but not irreducible there exists at least one $\fh$-invariant isotropic subspace $\hat E\subset V$. 
\begin{lm}\label{LE}   The following statements are true for any indecomposable subalgebra $\fh$ of $\fg_2^*\subset \fso(4,3)$.
Let $E\subset V$ be an $\fh$-invariant isotropic subspace.  
\begin{enumerate}
\item If $\dim  E=1$, then $\hat E(E):=\{v\in V\mid \forall\, e\in E: v\times e=0\}$ is a three-dimensional isotropic $\fh$-invariant subspace of $V$ containing $E$.
\item  If $\dim E=3$, then there exists a uniquely determined one-dimen\-sional isotropic $\fh$-invariant subspace $E_0\subset  E$ such that  $E=\hat E(E_0)$.
\item If $\dim E=2$ and if $b_1\times b_2\not=0$ for a basis $b_1,b_2$ of $E$, then there exists a one-dimensional $\fh$-invariant subspace $E_0\subset V$ not contained in $E$ such that $E\oplus E_0$ is isotropic.
\end{enumerate}
\end{lm}

\proof
(1) Suppose $E=\RR\cdot b$, $b\in V$. We want to show that $\hat E(E)=U(b\cdot\psi_0)$. Note first that $b\perp \hat E(E)$ since $0=b\times(v\times b)=\langle b,b\ra v-\la b.v\ra b=-\la b,v\ra b$ for all $v\in \hat E(E)$. Now (\ref{Ecross}) shows that $\hat E(E)\subset U(b\cdot \psi_0)$. Equation (\ref{Ecross}) also implies  
$\psi_0^\perp \ni (u\times b)\cdot\psi_0=\la u,b\ra\psi_0$ for all $u\in U(b\cdot \psi_0)$, thus  $U(b\cdot \psi_0)\subset \hat E(E)$. The subspace $U(b\cdot\psi_0)$ is three-dimensional and isotropic. Moreover, it is $\fh$-invariant since $E$ is $\fh$-invariant and $\psi_0$ is annihilated by $\fh$. Hence the same is true for $\hat E(E)$.

(2) Suppose $\dim E=3$. Then there exists an isotropic spinor $\ph_0$ such that $E=U(\ph_0)$. Since $\fh(E)\subset E$, we have $ E\cdot(\fh\cdot\ph_0)=\fh\cdot( E\cdot\ph_0)=0$, hence $E\cdot\ph_0\subset \RR\cdot \ph_0$. Assume that $\langle \ph_0,\psi_0\rangle_\Delta\not=0$. We define a vector $X$ by $\psi_0+X\cdot\psi_0\in\RR\cdot\ph_0$. Then $\fh(X)=0$ because of $\fh\cdot(\psi_0+X\cdot\psi_0)=\fh(X)\cdot\psi_0\in(\RR\cdot\ph_0)\cap \psi_0^\perp=0$. On the other hand, $X$ is not isotropic, which is a contradiction to indecomposability. Hence $\langle\psi_0,\ph_0\rangle_\Delta =0$. In particular, we can define $b\in V$ by $\ph_0=b\cdot\psi_0$. Then $E=U(\ph_0)=\hat E(E_0)$ for $E_0=\RR\cdot b$. As for uniqueness,  $\hat E(\RR\cdot b)=\hat E( \RR\cdot b')$ implies $U(b\cdot\psi_0)=U(b'\cdot\psi_0)$, thus $b'\in \RR\cdot b$.

(3) Now assume that $E=\Span\{b_1,b_2\}$ and that $b_1\times b_2=:b\not=0$. 
If $b$ were in $E$, then $b_1\times b_2=b_2$ without loss of generality. But this would imply $b_2=b_1\times(b_1\times b_2)=-\la b_1,b_1\ra b_2+\la b_1,b_2\ra b_1=0$, which contradicts $b_2\not=0$. Thus $E_0:=\RR\cdot b$ is not contained in $E$. Moreover, $E\oplus E_0$ is isotropic since $u\times v\perp u$ for all $u,v\in V$. 
\qed

Let $S$ be the socle of the holonomy representation. Then $S$ is isotropic. Indeed, $S\cap S^\perp$ is $\fh$-invariant. Since $\fh$ acts semisimply on $S$, there exists an $\fh$-invariant complement of $S\cap S^\perp$ in $S$. This complement is non-degenerate, hence trivial. Thus $S\cap S^\perp=S$.

\begin{de}
The holonomy representation is said to be of Type I, II or III if the dimension of $S$ equals one, two or three, respectively.
\end{de}

\subsection{Berger algebras of Type I}\label{SI}
Let $\fh$ be of Type I .
\begin{lm} \label{LI1} If $\fh$ is of type I, then there exists a basis $b_1,\dots,b_7$ of $V$ such that 
\begin{eqnarray*}
\ip&=&2(b^1\cdot b^5+b^2\cdot b^6+b^3\cdot b^7)- (b^4)^2\\
\omega&=&\sqrt 2(b^{167}+b^{235})-b^4\wedge (b^{15}-b^{26}-b^{37})
\end{eqnarray*} 
and $\fh$ is a subalgebra of
$$ \fh^I:=\{ h(A,v,u,y)\mid A\in\fgl(2,\RR),\ v\in\RR,\ u, y \in\RR^2\},$$
where 
$$h(A,v, u,y)= \left(
\begin{array}{ccccccc}
\tr A &-u_2&u_1&\sqrt 2 v&0&-y_1&-y_2\\
0&a_1&a_2&\sqrt 2 u_1&y_1&0&v\\
0&a_3&a_4&\sqrt 2 u_2&y_2&-v&0\\
0&0&0&0&\sqrt 2 v&\sqrt 2 u_1&\sqrt 2 u_2\\
0&0&0&0&-\tr A&0&0\\
0&0&0&0&u_2&-a_1&-a_3\\
0&0&0&0&-u_1&-a_2&-a_4
	\end{array}\right)$$ 
for	$ A=\left(\begin{array}{cc} a_1&a_2\\a_3&a_4\end{array}\right)$, $y=(y_1,y_2)^\top$, $u=(u_1,u_2)^\top$.
\end{lm}
\proof Let the socle be spanned by the isotropic vector $b$. Recall that $G_2^*$ acts transitively on isotropic lines in $\RR^{4,3}$. Hence we may assume $b=e_1$. Put $b_i:=e_i$, $i=1,\dots,7$. Then the assertion follows from (\ref{EG2}).
\qed

We define
$$\fm:=\{h(0,v,u,y)\mid v\in\RR,\ u, y \in\RR^2\}\subset \fh^I$$
and identify $\fgl(2,\RR)$ with $\{h(A,0,0,0)\mid A\in \fgl(2,\RR)\}$.
Then $$\fh^I=\fgl(2,\RR)\ltimes \fm.$$ 

We define the  matrices 
$$C_a:=\left(\begin{array}{cc} a&-1\\1&a\end{array}\right),\ S:=\left(\begin{array}{cc} 1&1\\0&1\end{array}\right),\ N:=\left(\begin{array}{cc} 0&1\\0&0\end{array}\right);
$$
and the following Lie subalgebras of $\fgl(2,\RR)\subset \fh^I$\,:
\begin{eqnarray*}
\fd&:=&\{\diag(a,d)\mid a,d\in\RR\},\\
\fu(1)&:=&\left\{ \left(\begin{array}{cc} a&-b\\b&a\end{array}\right) \ \Big| \ a,b\in\RR\right\},\\
\hat\fb_2&:=&\Span\{ I,N\},\\
\fs_\lambda&:=& \Span \{\diag(\lambda, \lambda-1),\ N\},\ \lambda\in\RR,\\
\fb_2&:=&\mbox{ Lie algebra of upper triangular matrices}.
\end{eqnarray*}

Furthermore, we define vector subspaces of $\fm$ by $\fm(0,0,0):=0$ and
\begin{eqnarray*}
\fm(1,0,0)&:=& \{h(0,v,0,0)\mid v\in\RR\},\\
\fm(0,1,0)&:=& \{h(0,0,(u_1,0)^\top, 0)\mid u_1\in\RR\},\\
\fm(0,0,1)&:=& \{h(0,0,0,(y_1,0)^\top)\mid y_1\in\RR\},\\
\fm(0,0,2)&:=& \{h(0,0,0, y)\mid y\in\RR^2\}.
\end{eqnarray*}
Now we put 
$$\fm(i,j,k)= \fm(i,0,0)\oplus\fm(0,j,0)\oplus \fm(0,0,k)$$
for $i,j\in\{0,1\}$, $k\in\{0,1,2\}$.

Let $\fa$ be the projection of $\fh$ to $\fgl(2,\RR)\subset \fh^I$. 
\begin{theo}\label{T1}
If $\fh$ is of Type I, then there exists a basis of $V$ such that we are in one of the following cases
\begin{enumerate}
\item $\fa\in \{0,\fsl(2,\RR),\, \fgl(2,\RR),\, \fu(1),\,\fb_2,\, \hat \fb_2,\, \fd,\,\RR\cdot C_a,\, \RR\cdot S\}$ and $\fh=\fa\ltimes \fm$,
\item $\fa=\fs_\lambda=\Span \{X:=\diag(\lambda,\lambda-1),\, N\}$ and  
\begin{enumerate}
\item $\lambda\in\RR$ and $\fh=\fa\ltimes \fm$, 
\item $\lambda=1$ and $\fh=\RR\cdot h(X,0,(0,1)^\top,0)\ltimes (\RR\cdot N\ltimes \fm(1,1,2)),$
\item $\lambda=2$ and $\fh=\Span\{X,\, h(N,0,(0,1)^\top,0) \}\ltimes\fm(i,j,2)$, where\\ $i,j\in\{0,1\},$
\end{enumerate}
\item $\fa=\RR\cdot\diag(1,\mu)$ and  
\begin{enumerate}
\item $\mu\in [-1,1]$ and $\fh=\fa\ltimes \fm$,
\item $\mu=0$ and $\fh=\RR\cdot h(\diag(1,0),0,(0,1)^\top,0)\ltimes \fm(1,1,2)$, 
\end{enumerate}
\item $\fa=\RR\cdot N$ and 
\begin{enumerate}
\item $\fh=\fa\ltimes \fm$, 
\item $\fh=\RR\cdot h(N,0,(0,1)^\top,0)\ltimes \fm(1,j,2)$ for $j\in\{0,1\}$.
\end{enumerate}
\end{enumerate}
\end{theo}

The remainder of this section is devoted to the proof of Theorem~\ref{T1}. Let us first have a closer look at the structure of $\fh^I=\fgl(2,\RR)\ltimes \fm$.
The element $A\in\fgl(2,\RR)$ acts on $\fm$ by 
$$
A\cdot h(0,v,u,y)=
h(0,\tr(A)v,Au,(A+\tr A)y).
$$
Furthermore, the Lie bracket on $\fm$ is given by
\begin{equation}\label{Elbm}
[h(0,v,u,y),h(0,\bar v,\bar u,\bar y)]=h(0,2\theta(u,\bar u),0,3(\bar v u-v\bar u)),
\end{equation}
where $\theta(u,\bar u):=u_1\bar u_2-u_2\bar u_1$
for $u,\bar u\in\RR^2$.

Similarly to $\fgl(2,\RR)$, we identify  $\GL(2,\RR)$ with a subgroup of $\G$ consisting  of block diagonal matrices: 
\begin{equation}\label{EGL}
\GL(2,\RR)\ni g\longmapsto \diag(\det g, g,1,(\det g)^{-1}, (g^\top)^{-1})\in \G,
\end{equation}
where $\G$ is considered with respect to the basis in Theorem \ref{T1}. Then
\begin{equation}\label{EAdA} 
\Ad(g)(h(A,v,u,y))= 
h(g A g^{-1},\, \det (g)\cdot v,\,gu,\,\det(g)\cdot gy).
\end{equation}
\begin{lm}\label{LaI}
Either $\fa\in\{0,\ \fsl(2,\RR),\ \fgl(2,\RR)\}$ or the basis $b_1,\dots, b_7$ in Lemma \ref{LI1} can be chosen such that $\fa$ is equal to one of the following Lie algebras:
\begin{enumerate}
\item $\RR\cdot A$, where $A$ is one of the matrices \\
$C_a,\, S,\, N,\, \diag(1,\mu),\ \mu\in [-1,1];$
\item $\fd$, $\fu(1)$, $\hat\fb_2$, $\fs_\lambda$, $\lambda\in\RR$;
\item $\fb_2$.
\end{enumerate}
\end{lm}
\proof  We identify $\GL(2,\RR)$ with a subgroup of $\G$ as described above. The conjugation of $\fh$  by an element of $\GL(2,\RR)$  is given by (\ref{EAdA}). Hence the proof of the Lemma is just the well-known classification of subalgebras of $\fgl(2,\RR)$ up to conjugation by $\GL(2,\RR)$:

(1) Suppose that $\fa$ is generated by one matrix $A\not=0$. If $A$ has real eigenvalues, then we may assume that one of the eigenvalues equals $1$. Otherwise, we may assume that the imaginary part of the eigenvalues equals $\pm 1$. Conjugating by an element of $\GL(2,\RR)$ we can achieve that $A$ has real Jordan normal form, which is $\diag(1,\mu),\,C_a,\,S$ or  $N$.

(2)  Now let $\fa$ be two-dimensional. If the natural representation of $\fa$ on $\RR^2$ is semisimple, then $\fa=\fd$ or $\fa=\fu(1)$ after conjugation by an element of $\GL(2,\RR)$. If not, then $\fa=\hat \fb_2$ or $\fa=\fs_\lambda$ after conjugation depending on whether $\fa$ is abelian or not.

(3) If $\dim \fa=3$ and $\fa\not=\fsl(2,\RR)$, then $\fa$ is solvable, thus conjugated to $\fb_2$.
\qed

We define
$${\mathcal K}(\fh)=\{R\in{\textstyle \bigwedge^2} V^*\otimes \fh \mid \forall x,y,z\in V: R(x,y)z+R(y,z)x+R(z,x)y=0\}$$
and 
$$\underline \fh:= \Span\{ R(x,y)\mid x,y\in V,\, R\in {\mathcal K}(\fh)\}.$$
Berger's first criterion implies $\fh=\underline\fh$.
Let $b_1,\dots,b_7$ be a basis as chosen in Lemma \ref{LI1}. If $R\in {\mathcal K}(\fh)$, then 
$$R_{ij}:=R(b_i,b_j)=h(A^{ij},v^{ij},u^{ij},y^{ij}).$$
\begin{lm} \label{Lhol1} If $R\in{\mathcal K}(\fh)$, then 
\begin{enumerate}
\item $R_{1j}=0$ for all $j\not=5$ and $R_{ij}=0$ for $i,j\in\{2,3,4\}$,
\item $\tr A^{ij}=0$ if $i<j$ and $(i,j)\not\in\{(5,6), (5,7)\}$,
\item $R_{15}=h(0,0,0, (\tr A^{56},\tr A^{57})^\top)$.
\end{enumerate}
\end{lm}
\proof 
Let $R$ be in ${\mathcal K} (\fh)$. Since $\la R_{ij}(b_k),b_l\ra=\la R_{kl}(b_i),b_j\ra$ and $R_{kl}\in\fh$, assertion
(1) follows. 

We define $b(i,j,k):=R_{ij}(b_k)+R_{jk}(b_i)+R_{ki}(b_j)$. From $b(i,j,5)=0$ we get $\tr A^{ij}=0$ for $i,j\not=5$. Furthermore, $b(1,5,6)=b(1,5,7)=0$ together with $R_{16}=R_{17}=0$ implies $A^{15}=0$, $u^{15}=0$, $v^{15}=0$ and $(\tr A^{56},\tr A^{57})^\top=y^{15}$. Now  $b(1,i,5)=0$ together with $u^{15}=0$ and $v^{15}=0$ gives $\tr A^{i5}=0$ for $i=2,3,4$. 
\qed

\begin{co} \label{Chol} If $\fa$ contains an element $A$ with $\tr A\not=0$, then 
$$\fy:= \{ y\in\RR^2 \mid h(0,0,0,y)\in\fh \}\not=0.$$
\end{co}
\begin{pr} \label{Pholtab1} The space ${\mathcal K}(\fh)$ can be parametrised by $a_i,\, r_i,\, x_i\in\RR$ ($i=1,2,3$), $b_k,c_k,, u_k,j_k\in\RR$ ($k=1,\dots,4$) and $v_1,v_2,t\in\RR$, where $R=h(A,v,u,y)\in {\mathcal K}(\fh)$ is given by the data in Table 1.
\begin{small}
\begin{table}
\renewcommand{\arraystretch}{1.4}
\begin{tabular}{|c|c|c|c|c|}
\hline
$R(b_i,b_j)$&$A$&$v$&$u$&$y$\\
\hline
$R_{15}$&$0$&$0$&$0$&$(b_1+b_4,c_1+c_4)$\\
$-R_{25}=\frac1{\sqrt2}R_{47}$&$\left(\begin{array}{cc} x_2&-x_1 \\x_3&- x_2\end{array}\right)$&$c_1-b_3$&$(r_2,r_3)$&$(u_2,u_4)$\\[3.5ex]
$R_{26}$&$\left(\begin{array}{cc} -a_1&-a_2 \\-a_3&a_1 \end{array}\right)$&$-r_2$&$(x_1,x_2)$&$(b_1,c_1)$\\[3.5ex]
$R_{27}$&$\left(\begin{array}{cc} -a_3& a_1\\j_1&a_3 \end{array}\right)$&$-r_3$&$(x_2,x_3)$&$(b_3,c_3)$\\[3.5ex]
$R_{35}=\frac1{\sqrt2}R_{46}$&$\left(\begin{array}{cc} x_1&x_4 \\x_2&-x_1 \end{array}\right)$&$b_4-c_2$&$(r_1,r_2)$&$(u_1,u_3)$\\[3.5ex]
$R_{36}$&$\left(\begin{array}{cc} -a_2& j_2\\a_1&a_2 \end{array}\right)$&$r_1$&$(x_4,-x_1)$&$(b_2,c_2)$\\[3.5ex]
$-R_{67}=\frac1{\sqrt2}R_{45}$&$\left(\begin{array}{cc} -r_2&r_1 \\-r_3&r_2 \end{array}\right)$&$u_2-u_3$&$(b_4-c_2,c_1-b_3)$&$(v_1,v_2)$\\[3.5ex]
$R_{56}$&$\left(\begin{array}{cc} b_1&b_2 \\b_3& b_4\end{array}\right)$&$v_1$&$(u_1,u_2)$&$(j_3,t)$\\[3.5ex]
$R_{57}$&$\left(\begin{array}{cc} c_1&c_2 \\c_3&c_4 \end{array}\right)$&$v_2$&$(u_3,u_4)$&$(t,j_4)$\\[3.5ex]
\hline
\multicolumn{5}{|c|}{$R_{12}=R_{13}=R_{14}=R_{16}=R_{17}=R_{23}=R_{24}=R_{34}=0$}\\[0.5ex]
\multicolumn{5}{|c|}{$R_{37}=R_{15}-R_{26}$}\\[0.5ex]
\hline
\end{tabular}
\\[0.5ex] 
\begin{center} 
    {{\bf  Table 1.}} 
\end{center}
\end{table}
\end{small}
\end{pr}
\proof Let $R$ be in ${\mathcal K} (\fh)$. As in the proof of Lemma~\ref{Lhol1}, we use $\la R_{ij}(b_k),b_l\ra=\la R_{kl}(b_i),b_j\ra$ and $R_{kl}\in\fh$, which now gives 
$$
{\sqrt2} R_{25}=-R_{47},\ {\sqrt2} R_{35}=R_{46},\ {\sqrt2}R_{67}=-R_{45}
$$
and  
$$
R_{37}=R_{15}-R_{26}.
$$
Let us consider the equations $b(i,j,k)=0$ for $i,j,k\not=4$. These equations give, in particular,
$$ v^{56} = -y_1^{67} =: v_1,\quad v^{56} = -y_1^{67} =: v_2,\quad y_1^{57} =y_2^{56}=:t.$$
Moreover, they imply the already proven properties of $R$ stated in Lemma \ref{Lhol1}.
The system of the remaining linear equations for the coefficients $A^{ij},\ v^{ij}$, $u^{ij}$ and $y^{ij}$ of $R_{ij}$ following from $b(i,j,k)=0$ for $i,j,k\not=4$ decomposes into five subsystems.  Each of these subsystems is a system of equations in the elements of one of the following sets:

\begin{align*}
M_1:= &\{ A^{26},A^{27}, A^{36}, A^{37}\},\\ 
M_2:= &\{A^{56},A^{57}, 
y^{26}, y^{36}, y^{27}, y^{37}, v^{25}, v^{35}, u^{67}\},\\
M_3:= &\{A^{25},A^{35}, u^{27}, u^{37}, u^{26}, u^{36}\},\\
M_4:= &\{A^{67}, u^{25}, u^{35}, v^{26}, v^{27}, v^{36}, v^{37}\},\\
M_5:= &\{u^{56}, u^{57}, y^{25}, y^{35}, v^{67}\}.
\end{align*}
The subsystem for $M_1$ is 
\begin{align*}
&-a_1^{26}=a_4^{26}=a_2^{27}=a_1^{37}=-a_4^{37}=:a_1,\quad a_4^{26}=a_3^{36},\\
&-a_2^{26}=-a_1^{36}=a_4^{36}=a_2^{37}=:a_2,\\
&-a_3^{26}=-a_1^{27}=a_4^{27}=a_3^{37}=:a_3.
\end{align*}
Together with $a_3^{27}=:j_1$, $a_2^{36}=:j_2$ this gives the parametrisation of $M_1$ claimed in the proposition.

For $M_2$, we have  
\begin{align*}
&a_1^{56}=y_1^{26}=:b_1, \quad a_2^{56}=y_1^{36}=:b_2,\\
&a_3^{57}=y_2^{27}=:c_3, \quad a_4^{57}=y_2^{37}=:c_4,\\
&v^{25}+y_2^{26}=a_3^{56}=a_1^{57}+u_2^{67}=:b_3,\\
&v^{35}+y_2^{36}=a_4^{56}=a_2^{57}-u_1^{67}=:b_4,\\
&v^{25}+a_1^{57}=y_1^{27}=y_2^{26}+u_2^{67},\\
&v^{35}+a_2^{57}=y_1^{37}=y_2^{36}-u_1^{67}
\end{align*}
and for $M_3$ 
\begin{align*}
&a_1^{35}=a_2^{25}=u_1^{26}=-u_2^{36}=-u_1^{37}=-a_4^{35}=:x_1,\\
&a_4^{25}=a_3^{25}=-u_2^{37}=u_1^{27}=u_2^{26}=-a_1^{25}=:x_2,\\
&u_2^{27}=-a_3^{25}=:x_3,\quad u_1^{36}=a_2^{35}=:x_4.
\end{align*}
Moreover, 
\begin{align*}
&a_1^{67}=-v^{26}=-u_1^{25}=u_2^{35}=v^{37}=-a_4^{67}=:r_2,\\
&-a_2^{67}=v^{36}=u_1^{35}=:r_1,\\
&a_3^{67}=-v^{27}=-u_2^{25}=:r_3
\end{align*}
is the system for $M_4$ and 
\begin{align*}
&y_1^{35}=u_1^{56}=:u_1, \quad -y_1^{25}=u_2^{56}=:u_2,\\
&y_2^{35}=u_1^{57}=:u_3, \quad -y_2^{25}=u_2^{57}=:u_4,\\
&v^{67}=u_1^{57}-u_2^{56}
\end{align*}
that for $M_5$. This shows that $R\in {\mathcal K}(\fh)$ is given as in Table 1. Moreover, the above considered systems of linear equations imply that each $R$ that is defined as in Table 1 for an arbitrary choice of the parameters $a_i, r_i, x_i,b_k,c_k, u_k,j_k,v_1,v_2,t$ satisfies $b(i,j,k)=0$ for $i,j,k\not=4$. 
It is easy to show by a direct calculation that it also satisfies $b(i,j,4)=0$ for $1\le i,j\le7$.
\qed

{\it Proof of Theorem \ref{T1}.} We define
\begin{eqnarray}
\fu&:=& \{u\in\RR^2 \mid \exists v\in\RR,\,\exists y\in\RR^2:\, h(0,v,u,y)\in \fh\}, \label{Edefu}\\
\fv&:=& \{v\in\RR \mid \exists  y\in\RR^2:\, h(0,v,0,y)\in \fh\} \label{Edefv}
\end{eqnarray}

Below, we will several times conjugate by elements of $\G$. In particular, we will use the following formulas, which can be obtained from (\ref{Elbm}): 
\begin{eqnarray}
\lefteqn{\Ad(\exp h(0,\bar v,0,0))( h(A,v,u,y))=h(A, v-\tr(A)\bar v, u, y-3\bar v u), } \label{EAdv} \\[1ex]
\lefteqn{\Ad(\exp h(0,0,\bar u,0))( h(A,v,u,y))= } \label{EAdu}\\
\nonumber
&& h(A, v-2\theta(u,\bar u)-\theta(\bar u, A\bar u), u-A\bar u, y+(3v -3\theta(u,\bar u)-\theta(\bar u, A\bar u))\bar u),\\[1ex]
\lefteqn{\Ad(\exp h(0,0,0,\bar y))(h(A,v,u,y))= 
h(A,v,u,y-(A+\tr A)\cdot \bar y).} \label{EAdy}
\end{eqnarray}
(1) If $\fa=0$, then $\dim \fu=2$ since otherwise the kernel of $\fh$ would be at least two-dimensional. Thus there exist elements $v_1,v_2\in \RR$,
$y_1,y_2\in\RR^2$ such that $$h_1:=h(0,v_1,(1,0)^\top,y_1),\ h_2:=h(0,v_2,(0,1)^\top,y_2)\in \fh.$$ Then also
$h_3:=[h_1,h_2]= h(0,2,0,3(v_2,-v_1)^\top)$, 
$[h_1,h_3]=h(0,0,0,(6,0)^\top)$ and
$[h_2,h_3]=h(0,0,0,(0,6)^\top)$
are elements of $\fh$, hence $\fh=\fm$.

Suppose $\fa=\fsl(2,\RR)$. Assume that $\fu=0$. Then for each $A\in\fa$ there exists a unique element $u_A\in\RR^2$ such that $h(A,v,u_A,y)\in\fh$ for some $v\in\RR, y\in\RR^2$. The map $A\mapsto u_A$ is a cocycle with respect to the standard representation of $\fa$ on $\RR^2$. Since $H^1(\fa,\RR^2)=0$ by Whitehead's lemma, it is a coboundary, i.e., $u_A=A\hat u$ for some $\hat u\in\RR^2$. After conjugation of $\fh$ by $\exp h(0,0,\hat u,0)$ according to (\ref{EAdu}), the projection of $\fh$ to $\{h(0,0,u,0)\mid u\in\RR^2\}$ is trivial. But then $\fh$ is of Type III. Hence $\fu\not=0$. Since $\fu$ is invariant under $\fa$, we obtain $\dim \fu=2$ and as above we conclude $\fm\subset\fh$. 

For $\fa=\RR\cdot C_a$ we can argue similarly. Indeed, $h(C_a,v,u,y)$ for some $v,u,y$. Since $C_a$ is non-singular, we have $u=C_a(\hat u)$ for some $\hat u\in\RR^2$ and we can proceed as above. 

Now suppose $\fa=\RR\cdot S$. Since $S$ defines a bijective map we may assume that (after conjugation) $h(S,v,0,y)\in\fh$ for some $v\in\RR$ and $y\in\RR^2$. Hence $\fu\not=0$ since otherwise $\fh$ would be of Type II. Assume that $\dim\fu=1$. Then $\fu=\RR\cdot (1,0)^\top$ since $\fu$ is invariant under $\fa$. But then $\fh$ again would be of Type II. Hence $\dim\fu=2$, which implies $\fm\subset\fh$.

Before we continue with the remaining cases let us make the following remark. Suppose $I\in\fa$ and $\fu=0$. Then, as above, $h(A,v,u_A,y)\in\fh$ for a cocycle $A\mapsto u_A$ and some $v\in\RR, y\in\RR^2$. Since $I\in\fa$, we have, in particular, $0=u_{[I,A]}=I(u_A)-A(u_I)$. Hence $u_A=A(u_I)$ for all $A\in \fa$. As in the case of $\fa=\fsl(2,\RR)$, after conjugation of $\fh$ by $\exp h(0,0,\hat u,0)$ the projection of $\fh$ to $\{h(0,0,u,0)\mid u\in\RR^2\}$ is trivial and $\fh$ is not of Type I. 

Suppose $\fa\in\{\fgl(2,\RR),\, \fu(1)\}$. Then $\fu\not=0$ by the above remark. Since $\fu$ is $\fa$-invariant we get $\dim \fu=2$, thus $\fm\subset \fh$.

Now let $\fa$ be one of the Lie algebras $\fb_2$, $\hat\fb_2$. Then $\fu\not=0$. If $\dim \fu=2$, then $\fm\subset\fh$. Thus we have to discuss the case  $\dim\fu=1$. Because of the $\fa$-invariance of $\fu$ we have $\fu=\RR\cdot (1,0)^\top$. Since $\tr I\not=0$ and $I$ and $I +\tr I$ act bijectively on $\RR^2$, we may assume that $h_1:=h(I,0,0,0)$ is in~$\fh$ (after conjugation according to (\ref{EAdv}) -- (\ref{EAdy})). Let $h(A,v_A,(0,u')^\top,y_A)$ be in $\fh$. Then also
$$
[h_1,h(A,v_A,(0,u')^\top,y_A)]= h(0,2v_A, (0,u')^\top, 3y_A)
$$
is in $\fh$, which implies $u'=0$. But then $\fh$ would be of Type II.  

Finally, we consider the case $\fa=\fd$. Then again $\fu\not=0$. If $\dim\fu=2$, we are done. Assume $\dim\fu=1$. Because of the invariance of $\fu$ under $\fa$ we have $\fu=\RR\cdot (1,0)^\top$ or $\fu=\RR\cdot (0,1)^\top$. Since $\tr I\not=0$ and $I$ and $I +\tr I$ acts bijectively on $\RR^2$, we may assume that $h_1:=h(I,0,0,0)$ is in~$\fh$ (after conjugation). For $H=\diag(1,-1)$, there are elements  $v\in\RR$, $u\in\RR^2$, $y\in\RR^2$ such that $h(H,v,u,y)\in\fh$. Since
$$[h_1,h(H,v,u,y)]=h(0,2v,u,3y)\in\fh,$$
$u\in\fu$ by (\ref{Edefu}). Thus we may assume that $h(H,v,0,y)\in\fh$ for some $v\in\RR$, $y\in\RR^2$. But then $\fh$ is of Type II. 

(2) Now let $\fa$ be equal to $\fs_\lambda$. First we want to show that $\fu=0$ implies $\lambda=2$. If $\fu=0$ we can define a cocycle $A\mapsto u_A$ by $h(A,v,u_A,y)\in\fh$ for some $v,y$. Since $[X,N]=N$ we obtain $X\cdot u_N-N\cdot u_X=u_N$, which yields
\begin{equation}\label{EXN}
(\lambda-1) u^1_N=u_X^2,\quad (\lambda-2) u_N^2=0
\end{equation}
for the components $u_X^1$, $u_X^2$ of $u_X$ and $u^1_N$, $u^2_N$ of $u_N$. We may assume that $u_X^2=0$. Indeed, if $\lambda\not=1$, then $\{(0,u_2)^\top\mid u_2\in\RR\}$ is in the image of $X$, thus we can find a suitable conjugation of $\fh$. If $\lambda=1$, then $u_X^2=0$ follows from (\ref{EXN}). If $\lambda\not=2$, then (\ref{EXN}) gives $u_N^2=0$ and, consequently, the projection of $\fh$ to $\{h(0,0,(0,u_2)^\top,0)\mid u_2\in\RR\}$ vanishes. Hence $\fh$ is of Type~II. Thus $\fu=0$ can hold only if $\lambda=2$. 

Now take, $\lambda\in\RR$, $\lambda\not\in\{1,2\}$. We have already seen that $\fu\not=0$. If $\dim \fu=2$, we are done. Assume that $\dim\fu=1$. Then $\fu=\RR\cdot (1,0)^\top$ by $\fa$-invariance of $\fu$. Since $\lambda\not=1$, we see as above that we may assume that $h(X,v_X,u_X,y_X)$ is in $\fh$ for some $v_X,u_X,y_X$ with $u_X=(u',0)$. Hence $h_X:=h(X,v_X,0,y_X)\in\fh$. Furthermore, $h_N:=h(N,v_N,u_N,y_N)\in \fh$ for some $v_N,u_N,y_N$ with $u_N=(0,u'')^\top$. Then 
$$[h_X,h_N]=h(N,(2\lambda-1)v_N,(\lambda-1)u_N, y)$$ for some $y\in\RR^2$. Since $\lambda\not=2$, we obtain $u''=0$. Hence the projection of $\fh$ to $\{h(0,0,(0,u_2)^\top,0)\mid u_2\in\RR\}$ vanishes and $\fh$ is of Type II, which contradicts our assumption.

Suppose now $\lambda=1$. We already know that $\fu\not=0$. If $\dim \fu=2$, then we are in case (2)(a). Let us consider the case $\dim\fu=1$. Then $\fu=\RR\cdot (1,0)^\top$. Choose $v_X,u_X,y_X,v_N,u_N,y_N$ such that $h_X:=h(X,v_X,u_X,y_X)$ and $h_N:=(N,v_N,u_N,y_N)$ are in $\fh$. Since $\tr X\not=0$ and $X+\tr X$ acts bijectively, we may assume $v_X=0$ and $y_X=0$. Then $[h_X,h_N]=h(N,v_N,Xu_N,(X+\tr X)y_N)$, which implies $u_N-Xu_N\in\fu$. Since also $X u_N$ is in $\fu$, we see that $u_N\in\fu$. Thus we can choose $u_N=0$. To summarise, we get 
$$h_X=h(X,0, (0,u')^\top, 0),\quad h_N=h(N,v_N,0,y_N),$$
where $u'\not=0$ since otherwise $\fh$ would not be of Type I. We choose $v_0,y_0$ such that $h_0:=h(0,v_0,(1,0)^\top,y_0)\in\fh$. Then $[h_X,h_0]=h(0,v_0-2u', (1,0)^\top,y)$ for some $y\in\RR^2$. Since $u'\not=0$, this implies $\fv\not=0$. Hence, there exists $\hat y=(\hat y_1,\hat y_2)\in\RR^2$ such that $h_v:=h(0,1,0,\hat y)\in\fh$. We have 
$$[h_X,h_v]=h(0,1,0,(2\hat y_1, \hat y_2+3u')^\top).$$
Hence $(\hat y_1, 3u')\in\fy$. Since $\fy$ is $\fa$-invariant and $u'\not=0$ we obtain $\fy=\RR^2$. Consequently, $\fh=\RR\cdot h_X\ltimes (\RR\cdot N\ltimes \fm(1,1,2))$. Conjugating by $\diag (u', (u')^{-1})\in\SL(2,\RR)\subset \exp \fh\subset \SO(4,3)$  we get $u'=1$ and we are in case (2)(b).

Finally, suppose $\lambda=2$. Since $X=\diag(2,1)$, $\tr X\not=0$ and $X+\tr X$ acts bijectively, we have $X\in\fh$ after a suitable conjugation. Cor.~\ref{Chol} together with the $\fa$-invariance of $\fy$ implies that $\fy$ contains $\RR\cdot (1,0)^\top$. In particular, $h_N:=h(N,v_N,u_N, y_N)\in\fh$ for some $v_N, u_N=(u_N^1,u_N^2)^\top$ and $y_N=(0,y')^\top$. Then
\begin{equation}\label{EXhN}
[X,h_N]=h(N,3v_N, (2u_N^1,u_N^2)^\top, (0,4y')^\top).
\end{equation}
Hence $(u_N^1,0)\in\fu$, thus we can choose $u_N^1=0$. Let us first consider the case $\fv=0$. Then $h_N=h(N,0,(0,u')^\top, (0,y')^\top)$ and $h_v=h(0,1,0,y_v)$ are in $\fh$ for some $y_v\in\RR^2$, thus also $[h_N,h_v]=h(0,0,0,(y',3u')^\top)\in\fh$. If now $u'\not=0$, this shows $\fy=\RR^2$. Conjugation by $(u')^{-1}\cdot I\in\GL(2,\RR)\subset \SO(4,3)$ shows that we may assume $u'=1$. Hence we are in case (2)(c) with $i=1$ or  $\fh=\fa\ltimes \fm$. If $u'=0$, then $\fu=\RR^2$ since $\fh$ is of Type I and $\fu$ is $\fa$-invariant. Hence $\fh=\fa\ltimes \fm$. If $\fv=0$, then (\ref{EXhN}) implies $v_N=0$.  Assume $\fy=\RR\cdot (1,0)^\top$. Then (\ref{EXhN}) would imply $y'=0$. Hence the projection of $\fh$ to $\{ h(0,v,0,(0,y_2)^\top)\mid v,y_2\in\RR\}$ would be trivial. Thus Prop.~\ref{Pholtab1} would imply $c_1=\ldots=c_4=0$ and $b_4=0$. But then $X\not\in\underline \fh$, which would mean that $\fh$ is not a Berger algebra. Hence $\fy=\RR^2$. If $u'=0$, then $\fu=\RR^2$ since $\fh$ is of Type I and $\fu$ is $\fa$-invariant. Thus $\fh=\fa\ltimes \fm$. If $u'\not=0$ we again may assume $u'=1$. Then $\fh=\fa\ltimes \fm$ or we are in case (2)(c) with $i=0$.

(3) Let $\fa$ be spanned by $\diag(1,\mu)$. 
Assume first that $\mu\not=0$. Then, possibly after a further conjugation, $h(\diag(1,\mu),\hat v,0,\hat y)\in\fh$ for some $\hat v\in\RR$, $\hat y\in\RR^2$. If $\mu\not=1$, then $\RR\cdot (1,0)^\top$ and $\RR\cdot (0,1)^\top$ are the only proper $\fa$-invariant subspaces of $\fu$. Since $\fh$ is of Type I, we obtain $\fu=\RR^2$, thus $\fh=\fa\ltimes\fm$. If $\mu=1$, then the operators $h(0,v,u,y)$ with $u\in\fu$ do not have a non-trivial common kernel since $\fh$ is of Type I. Hence $\dim\fu=2$, which implies $\fh=\fa\ltimes\fm$. Now we consider $\mu=0$. Then we may assume that $\hat h:=h(\diag(1,\mu),0,(0,u')^\top,0)\in\fh$ for some $u'\in\RR$. Since $\fh$ is of Type I, we have $\RR\cdot (1,0)^\top\subset \fu$.  If $\fu=\RR^2$, then $\fh=\fa\times\fm$. If $\fu=\RR\cdot(1,0)^\top$, then $u'\not=0$ since $\fh$ is of Type I. Thus we may assume $u'=1$. Furthermore, $h_0:=h(0,v_0,(1,0)^\top,y_0)\in\fh$ for some $v,y_0$. Since $[\hat h,h_0]=h(0,v_0-2,(1,0)^\top,y)$ for some $y\in\RR^2$, we see that $\fv\not=0$. Hence, $\hat h= h(\diag(1,\mu),0,(0,1)^\top,0),\, h_0=h(0,0,(1,0)^\top,y_0),\, h_v:=h(0,1,0,y_v)$ are in $\fh$ for some $y_0,y_v=(y_v^1,y_v^2)^\top\in\RR^2$. Since $[\hat h, h_v]=h(0,1,0,(2y_v^1,y_v^2)^\top)$, we get $\RR\cdot (0,1)^\top \subset \fy$ and $[h_0,h_v]=h(0,0,0,(3,0)^\top)$ implies $\RR\cdot (1,0)^\top \subset \fy$. Thus $\fy=\RR^2$ and we are in case (3)(b).

(4) Finally, we consider $\fa=\RR\cdot N$. If $\fu=\RR^2$, then $\fh=\fa\ltimes\fm$. Suppose that $\fu\not=\RR^2$. Then $\fu=0$ or $\fu=\RR\cdot (1,0)^\top$.  Possibly after conjugation, $h_N:=h(N,v_N,(0,u')^\top, (0,y')^\top)\in\fh$ for suitable $v_N,u',y'\in\RR$. Then $u'\not=0$ since $\fh$ is of Type I. We may assume $u'=1$. Then, again after conjugation, $v_N=y'=0$. Let us first consider the case $\fv\not=0$. Then $h_v:=h(0,1,0,y_v)\in\fh$ for some $y_v\in\RR^2$. Because of $[h_N,h_v]=[h(N,0,(0,1)^\top,0),h_v]=h(0,0,0,(y_1,3)^\top)$ for some $y_1\in\RR$. Since $\fy$ is $\fa$-invariant, we obtain $\fy=\RR^2$. Hence $\fh$ is as claimed in case (4)(b). Now suppose $\fv=0$. If $h_0:=h(0,v_0,(1,0)^\top,y_0)$ would be in $\fh$ for some $v_0,y_0$, then also $[h_N,h_0]=[h(N,0,(0,1)^\top,0),h_0]=h(0,-2,0,\hat y)\in\fh$ for some $\hat y\in\RR^2$, which would contradict $\fv=0$. Thus $\fu=\fv=0$. This implies that all parameters describing ${\mathcal K}({\fh})$ are zero except of $j_3,j_4,t$, see Prop.~\ref{Pholtab1}. But then $h_N=h(N,0,(0,1)^\top,0)$ would not be in $\underline \fh$, which would contradict Berger's first criterion.
\qed

\subsection{Berger algebras of Type III} 
\begin{lm} If $\fh$ is of Type III, then there exists a basis $b_1,\dots,b_7$ of $V$ such that the metric on $V$ equals $2\sigma^1\cdot\sigma^5+2\sigma^2\cdot \sigma^6+2\sigma^3\cdot \sigma^7- (\sigma^4)^2$ with respect to the dual basis and $\fh$ is a subalgebra of
$$\fh^s = \{ h^s(A,v,y):=h(A,v,0,y)\mid A\in\fgl(2,\RR),\ v\in\RR,\ y \in\RR^2\}\subset \fh^I.$$
\end{lm}
Let $\fa$ and $\fm(i,j,k)$ be defined as in Section~\ref{SI}.
\begin{theo}\label{T3}
If $\fh$ is of Type III, then there exists a basis such that we are in one of the following cases
\begin{enumerate}
\item $\fa\in \{\fsl(2,\RR),\, \fgl(2,\RR),\, \fu(1),\fd\}$ and $\fh=\fa\ltimes \fm(1,0,2)$ 
\item $\fa\in \{0,\, \RR\cdot\diag(1,0)\}$ and $\fh=\fa\ltimes \fm(1,0,k)$ for $k\in\{1,2\}$.
\end{enumerate}
\end{theo}
\proof If the socle $S$ is three-dimensional, it defines a one-dimensional invariant subspace in a natural way, see Lemma~\ref{LE}. Since $\G$ acts transitively on isotropic lines we may assume that this space is spanned by $e_1$. Then $S=\Span\{e_1,e_2,e_3\}$, see Lemma~\ref{LE}. Now we take the same basis as in Lemma~\ref{LI1}., i.e., $b_i=e_i$, $i=1,\dots,7$. Then $\fh\subset\fh^I$. Since the representation of $\fh$ on $S$ is semisimple,  $\RR\cdot b_1$ has an invariant complement $\hat S$ in $S$. Hence $\fu=\{u\in\RR^2\mid \exists v\in\RR,\,\exists y\in\RR^2:\, h(0,v,u,y)\in \fh\} =0$. Since $\fh$ acts semisimply on $\hat S$, $\fa\cong\GL(2,\RR)$ acts semisimply on $\RR^2$. Thus $\fa\in\{0,\fgl(2,\RR),\fsl(2,\RR)\}$ or $\fa$ is conjugated to one of the Lie algebras $\fu(1), \fd, \RR\cdot\diag(1,\mu)$ or $\RR\cdot C_a$. So we may assume that $\fa$ is one of these Lie algebras. In the proof of Theorem~\ref{T1}, we have seen that if $\fa\in\{\fsl(2,\RR),\RR\cdot C_a\}$ or if $I\in\fa$, then $\fu=0$ implies that, after a suitable conjugation, the projection of $\fh$ to $\{h(0,0,u,0)\mid u\in\RR^2\}$ is trivial, which means that $\fh\subset\fh^s$. Furthermore, if $\fa=\RR\cdot\diag(1,\mu)$, $\mu\not=0$, then we also have $\fh\subset\fh^s$ after a suitable conjugation according to (\ref{EAdu}) since $\diag(1,\mu)$ acts bijectively. If $\mu=0$, then $h(\diag(1,\mu),v_0,(0,u'),y_0)\in\fh$ after conjugation according to (\ref{EAdu}) since $\RR\cdot (1,0)^\top$ is in the image of $\diag(1,0)$. Since $\fh$ acts semisimply on $S=\Span\{b_1,b_2,b_3\}$, it follows that $u'=0$. 

Now we are using Berger's criterion in order to determine $\fh$.
Since $\fh\subset\fh^s$, Prop.~\ref{Pholtab1} implies 
$0=r_1=r_2=r_3=x_1=\dots=x_4=u_1=\dots=u_4$ and $b_4=c_2$, $c_1=b_3$ for the parameters of ${\mathcal K}(\fh)$. 
Since $\fh$ is indecomposable, we have $\fv\not=0$ since otherwise the non-isotropic vector $b_4$ would be in the kernel of $\fh$. Hence $v_1\not=0$ or $v_2\not=0$, which implies $\fy\not=0$. For $\fa\in\{\fgl(1,\RR), \fsl(2,\RR), \fu(1)\}$ we get $\fy=\RR^2$ since there is no one-dimensional invariant subspace of $\RR^2$. Thus $\fa\subset\fh$ and $\fh\cap\fm=\fm(1,0,2)$. If $\fa\subset\fd$, then the parameters of ${\mathcal K}(\fh)$ satisfy in addition  $0=a_1=a_2=a_3=j_1=j_2=0$ and $b_2=b_3=b_4=c_1=c_2=c_3=0$. Hence all parameters appearing in some $A^{ij}$ in Table 1 are zero except of $b_1$ and $c_4$. For $\fa=\fd$ this immediately implies $\fy=\RR^2$. For $\fa=\RR\cdot\diag (1,\mu)$ we obtain $\mu=0$ and $b_1\not=0$. Since $\diag(1,0)+ \tr \diag(1,0)$ act bijectively on $\RR^2$ we may have $\fa\subset\fh$ after conjugation according to (\ref{EAdy}). Furthermore, $b_1\not=0$ yields $\fh\cap\fm=\fm(1,0,k)$ for $k=1,2$, see Table 1. If $\fa=0$, then $\fy=\RR^2$ or $\fy=\RR\cdot (1,0)^\top$ after conjugation, which proves the assertion for this case. 

It remains to exclude $\fa=\RR\cdot C_a$.  For $\fa=\RR\cdot C_a$, Prop.~\ref{Pholtab1} implies the system 
$$b_1=b_4=c_2=-c_3,\ -b_2=b_3=c_1=c_4,\ b_1=ab_3,\ c_1=ac_3$$ 
of linear equations, which has only the trivial solution. This gives a contradiction to $h(C_a,0,0,y_0)\in\underline \fh$ for some $y_0\in\RR^2$.

\qed

\subsection{Berger algebras of Type II} Let $\fh$ be of Type II.

For $z=(z_1,\dots,z_4)\in\RR^4$, we define 
$$\sigma(z):= \left(\begin{array}{ccc} z_2&\sqrt 2 z_3&z_4\\z_1&\sqrt 2 z_2&z_3\end{array}\right), \qquad  \sigma(z)^*= \left(\begin{array}{cc} -z_4& -z_3\\\sqrt 2 z_3 &\sqrt 2 z_2\\ -z_2&-z_1\end{array}\right) $$
and, for $ A=\left(\begin{array}{cc} a_1&a_2\\a_3&a_4\end{array}\right)\in\fgl(2,\RR)$ and $c\in\RR$, we put
$$\rho(A):=\left(\begin{array}{ccc} a_1-a_4&-\sqrt 2 a_2&0\\-\sqrt 2 a_3&0&-\sqrt 2 a_2\\
0&-\sqrt 2 a_3& -a_1+a_4\end{array}\right),\qquad U(c):=\left(\begin{array}{cc} 0&-c\\c&0\end{array}\right).$$

\begin{lm} \label{LII} If $\fh$ is of type II, then there exists a basis $b_1,\dots,b_7$ of $V$ such that 
\begin{eqnarray*}
\ip&=&2(b^1\cdot b^6+b^2\cdot b^7+b^3\cdot b^5)- (b^4)^2,\\
\omega&=& \sqrt 2(-b^{157}+b^{236})-b^4\wedge(b^{16}-b^{27}-b^{35})
\end{eqnarray*}
 and $\fh$ is a subalgebra of
$$ \fh^{II}:=\{ h(A,z,c)\mid A\in\fgl(2,\RR),\ z \in\RR^4,\ c\in\RR\},$$
where
$$h(A,z,c)= \left(
\begin{array}{ccc}
A&\sigma(z)&U(c)\\
0&\rho(A)&\sigma(z)^*\\
0&0&-A^\top
\end{array}\right).$$ 
\end{lm}
\proof Let $b_1,b_2$ be a basis of the socle $S$. Then $b_1\times b_2=0$ since otherwise $S$ would be contained in a 3-dimensional isotropic subspace on which $\fh$ acts semisimply, see Lemma \ref{LE} (3). Since $\G$ acts transitively on isotropic vectors, we may assume $b_1=e_1$.  Because of $b_1\times b_2=0$, the vector $b_2$ is in $\Span\{e_1,e_2,e_3\}$. Furthermore, the subgroup $\GL(2,\RR)\subset \G$ defined by (\ref{EGL}) acts as $\GL(2,\RR)$ on $\Span\{e_2,e_3\}$, thus we may assume $b_2=e_2$. Moreover, we put $b_3:=e_3$, $b_4:=e_4$, $b_5:=e_7$, $b_6:=e_5$, $b_7:=e_6$.
\qed

Note that $\omega$ and $\ip$ with respect to the chosen basis differ from those that we considered in the section on Type I.
We will also consider another embedding of $\fgl(2,\RR)$ into $\fg_2^*$. In this subsection we identify $\fgl(2,\RR)$ with 
\begin{equation}\label{EglII}
\{h(A,0,0)\mid A\in \fgl(2,\RR)\}\cong\fgl(2,\RR)
\end{equation}
and define $\fa$ to be the projection of $\fh$ to $\fgl(2,\RR)\subset\fh^{II}.$ 
We set
$$\fn:=\{h(0,z,c)\mid z\in\RR^4,\, c\in\RR\}$$
and, for $i,j,k\in\{1,2,3,4\}$, we define
\begin{eqnarray*}
\fn(i,j)&:=&\{h(0,z,c)\mid z\in\RR^4,\, z_l=0 \mbox{ if } l\not\in\{i,j\},\ c\in\RR\}\,,\\
\fn(i,j,k)&:=&\{h(0,z,c)\mid z\in\RR^4,\, z_l=0 \mbox{ if } l\not\in\{i,j,k\},\ c\in\RR\}\,.
\end{eqnarray*}
\begin{theo}\label{T2}
If $\fh$ is of Type II, then there exists a basis of $V$ such that we are in one of the following cases
\begin{enumerate}
\item $\fa\in \{\fsl(2,\RR),\, \fgl(2,\RR)\}$ and $\fh=\fa\ltimes \fn$,
\item $\fa\in \{\fu(1),\, \RR\cdot C_a\}$ and $\fh=\fa\ltimes \fn$ or 
\begin{center}$\fh=\fa\ltimes \{ h(0, (3r, s, r, 3s), c)\mid r,s,c\in\RR\}.$\end{center}
\item $\fa=\fd$ and $\fh=\fa\ltimes\fn_1$, where 
\begin{center}$\fn_1\in\{ \fn,\, \fn(1,3),\, \fn(2,3),\, \fn(1,2,3),\, \fn(1,2,4)\}$,\end{center}
\item $\fa=\RR\cdot\diag(1,\mu)$, $\mu\in[-1,1)$, and
\begin{enumerate} 
\item $\mu\in[-1,1)$  and $\fh=\fa\ltimes \fn_1$, where
\begin{center} $\fn_1\in\{ \fn,\, \fn(2,3),\,\fn(1,2,3),\,\fn(1,2,4),\,\fn(1,3,4),\, \fn(2,3,4)\}$, \end{center} 
\item $\mu=1/2$ and $\fh=\RR\cdot h(\diag(1,1/2), (1,0,0,0),0)\ltimes \fn_1$, \\
where 
$\fn_1\in\{\fn(2,3), \fn(2,3,4)\}$,
\item $\mu=0$ and $\fh=\fa\ltimes\fn(2,4)$ or 
\begin{center}
$\fh=\RR\cdot h(\diag(1,0), (0,1,0,0),0)\ltimes \fn_1$,
\end{center}
 where 
$\fn_1\in\{\fn(1,4),\,\fn(3,4),\, \fn(1,3,4)\}$,
\end{enumerate}

\item $\fa\in\{0,\,\RR\cdot I\}$ and $\fh=\fa\ltimes \fn_1$, where
\begin{center}$\fn_1\in\{ \fn,\, \fn(1,3),\, \fn(2,3),\, \fn(1,3,4),\, \fn(2,3,4)\}$,\end{center}
or  $\fn_1$ is one of the Lie algebras $\{h(0,z,c)\mid z\in Z,\, c\in\RR\} $ for 
\begin{enumerate} 
\item \label{Ia} $Z=\{(z_1,0, z_1,z_4)\mid z_1,z_4\in\RR\}$,
\item \label{Ib}  $Z=\{(0, z_2,z_3,-z_2)\mid z_1,z_4\in\RR\}$,
\item \label{Ic} $Z=\{(z_1,\alpha z_1,\alpha z_4,z_4)\mid z_1,z_4\in\RR\}$, $\alpha\in\left[\frac{\sqrt 3-1}{\sqrt6},\frac{\sqrt3+1}{\sqrt 6}\right]$,  
\item \label{Id} $Z=\{(s z_1,\alpha z_2 ,-\alpha z_1,-z_2)\mid z_1,z_2\in\RR\}$, $s\in(0,1]$, $\alpha\in\RR$ such that  $3\alpha^2-(s+1)\alpha-s=0$,
\item  \label{Ie} $Z=\{(z_1,z_2,\kappa z_1,z_4)\mid z_1,z_2,z_4\in\RR\}$, $\kappa=\pm1$.
\end{enumerate}
\end{enumerate}
\end{theo}

The remainder of this section is concerned with the proof of Theorem~\ref{T2}. Let us first describe the structure of $\fh^{II}=\fgl(2,\RR)\ltimes \fn$.
The Lie bracket on $\fn$ is given by
$$[h(0,z,c),h(0,\hat z,\hat c)]=h(0,0,\eta(z,\hat z)),$$
where
$$\eta(z,\hat z)=-z_1\hat z_4+z_4\hat z_1+3z_2\hat z_3-3z_3\hat z_2.$$
Moreover, $\fh^{II}=\fgl(2,\RR)\ltimes \fn$, where $A\in\fgl(2,\RR)$ acts on $\fn$ by 
$$
A\cdot h(0,z,c)= h(0,A\cdot z, \tr(A)\cdot c),
$$
where the representation of $\fgl(2,\RR)$ on $\RR^4$ is given by the equation
$$\sigma(A\cdot z)=A\circ \sigma(z)- \sigma(z)\circ \rho(A).$$ In particular, the basis vectors $I\in\fgl(2,\RR)$ and 
$$H=\left(\begin{array}{cc} 1&0\\0&-1\end{array}\right),\ X=\left(\begin{array}{cc} 0&1\\0&0\end{array}\right),\ Y=\left(\begin{array}{cc} 0&0\\1&0\end{array}\right)$$
act by $I\cdot z=z$, 
$H\cdot z=(-3z_1,-z_2,z_3,3z_4)$, $X\cdot z=(0,z_1,2z_2,3z_3)$ and $Y\cdot z=(3z_2,2z_3,z_4,0).$ This representation integrates to a representation of $\GL^+(2,\RR)$. Putting $\diag(1,-1)\cdot z=(z_1,-z_2,z_3,-z_4)$, we finally get a representation of $\GL(2,\RR)$ on $\RR^4$, which we will consider in the following.

Let $P_d$ denote the space of homogeneous polynomials of degree $d$ in $x,y$. On $P_d$, we consider the representation of $\GL(2,\RR)$ given by $(A\cdot p)(x,y):=p((x,y)A)$. Then 
\begin{eqnarray*}
\phi_1:\ \RR^4 & \longrightarrow & P_3\\
 z=(z_1,\dots,z_4)& \longmapsto & z_1y^3+3z_2xy^2+3z_3x^2y+z_4x^3
 \end{eqnarray*}
is an isomorphism satisfying
$$\phi_1(A\cdot z)= \lambda_1^{-1}(A)\cdot \phi_1(z),$$
where $\lambda_1$ denotes the automorphism $$\textstyle{\lambda_1:\GL(2,\RR)\longrightarrow\GL(2,\RR),\quad A\longmapsto\det(A) A}.$$
In particular, $\phi_1$ maps $\GL(2,\RR)$-orbits in $\RR^4$ to $\GL(2,\RR)$-orbits in $P_3$.

The representation of $\GL(2,\RR)$ on $\RR^4$ that we considered above induces representations of $\GL(2,\RR)$ on $\bigwedge^n\RR^4$. We consider these representations for $n=2,3$. Let $e_1,\dots,e_4$ be the standard basis of $\RR^4$. Let us start with $n=2$. The complementary subspaces $W_0:=\Span\{w_0:=e_{23}-3e_{14}\}$ and $W':=\Span\{e_{12}, e_{13}, w':=e_{23}+3e_{14}, e_{24}, e_{34}\}$ of $\bigwedge^2\RR^4$ are invariant under $\GL(2,\RR)$. The isomorphism $\phi_2: W'\rightarrow P_4$ defined by 
\begin{eqnarray*}
\phi_2:\quad e_{12}&\longmapsto & y^4\\
e_{13}&\longmapsto & 2xy^3\\
w'&\longmapsto & 6x^2y^2\\
e_{24}&\longmapsto & 2x^3y\\
e_{34}&\longmapsto & x^4
\end{eqnarray*}
satisfies
$$\phi_2(A\cdot u)= \lambda_2^{-1}(A)\cdot \phi_2(u),$$
where $\lambda_2$ denotes the automorphism $$\textstyle{\lambda_2:\GL(2,\RR)\longrightarrow\GL(2,\RR),\quad A\longmapsto \sgn(\det A)\sqrt{|\det(A)|} A}.$$

For $n=3$ the situation is even simpler. The representation  $\bigwedge^3\RR^4$  is equivalent to $P_3$. An equivalence is given by
\begin{equation}\label{Eaequ}
\Phi_3: \ e_{123}\longmapsto y^3,\ \ e_{124}\longmapsto xy^2,\ \ e_{134}\longmapsto x^2y,\ \ e_{234}\mapsto x^3. 
\end{equation}

Next we determine the orbits of the $\GL(2,\RR)$-action on the projective spaces $\PP(P_3)$ and $\PP(P_4)$. The line spanned by a polynomial $p$ is denoted by $[p]$.
\begin{lm}\label{Lpol}
The elements $[x^3]$, $[x^2y]$ and $[x(x^2\pm y^2)]$ constitute a complete system of representatives of the orbit space $\GL(2,\RR)\setminus \PP(P_3)$. 

The following elements constitute a complete system of representatives of $\GL(2,\RR)\setminus \PP(P_4)$: $[x^4]$, $[x^3y]$, $[x^2y^2]$, $[xy(x^2+rxy+y^2)]$ for $r\in[0,\,3/\sqrt 2]$, and $[(x^2+y^2)(x^2+sy^2)]$ for $s\in[0,1]$.
 \end{lm}
\proof The first assertion is well known and easy to prove. Let us check the second one. We denote by $\la p \ra$ the orbit of $[p]$. By a zero of $p=p(x,y)$ we mean a (real) zero of $p$ on $\RR P^1$. If $p$ has a zero of multiplicity three or four, then $\la p\ra=\la x^4\ra$ or $\la p\ra=\la x^3y\ra$. If $p$ has two zeroes of multiplicity two, then $\la p\ra=\la x^2y^2\ra $. Let us now consider the remaining cases.

(1) Suppose that $p$ has exactly two or four zeroes and all zeroes are simple. Then
$$\la p\ra=\la xy(x^2+bxy +cy^2)\ra=\la xy(x^2+rxy +\kappa y^2)\ra, \quad \kappa=\pm1,$$
by rescaling $x$ and $y$. We want to show that we may assume $\kappa=1$. If $\kappa=-1$, then 
$x^2+rxy +\kappa y^2=x^2+rxy - y^2$  has two zeroes, hence
$$x^2+rxy - y^2=(x+qy)(x-\textstyle{\frac1q} y),\quad |q|\ge1.$$
We choose $\hat q\in\RR$ such that $\hat q^2 =q^2 +1$ and put $x=\hat x+\hat q \hat y$, $y=q\hat x$. Then 
\begin{eqnarray*}
\lefteqn{[xy(x+qy)(x-y/q)]=[(\hat x +\hat q\hat y)\hat x (\hat x +\hat q\hat y+q^2\hat x)\hat y]}\\
&&=[(\hat x +\hat q\hat y)\hat x (\hat q^2\hat x +\hat q\hat y)\hat y]
\, =\, [\hat x \hat y(\hat x +\hat q\hat y)(\hat x+\hat y/\hat q)]\\
&&=[ \hat x\hat y (\hat x^2+\hat r \hat x \hat y +\hat y^2)]\,. 
\end{eqnarray*}
Thus $\la p\ra=\la xy(x^2+rxy +y^2)\ra$, where $r\ge 0$ (otherwise replace $x$ by $-x$) and $r\not=2$. We distinguish two cases.

(1)(a) If $r\in[0,2)$, then $p$ has exactly two zeroes. Polynomials with different values for $r$ belong to different $\GL(2,\RR)$-orbits since any transformation that maps $[xy(x^2+rxy+y^2)]$ to $[xy(x^2+\hat rxy+y^2)]$ leaves invariant the set $\{[1:0],\,[0:1]\}\subset \RR P^1$ of zeroes. 

(1)(b) If $r>2$, then $p$ has four different zeroes. Thus
$\la p \ra=\la xy (x+qy)(x+y/q)]$, where $0<q<1$. If $q<1/\sqrt 2$, then we choose $\hat q\in (1/\sqrt 2,1)$ such that $q^2+\hat q^2=1$. For $x=\hat x +\hat q\hat y$ and $y=-q\hat x$, we get
\begin{eqnarray*}
\lefteqn{[xy(x+qy)(x+y/q)]=[(\hat x +\hat q\hat y)\hat x (\hat x +\hat q\hat y-q^2\hat x)\hat y]}\\
&&=[(\hat x +\hat q\hat y)\hat x (\hat q^2\hat x +\hat q\hat y)\hat y]
\, =\, [\hat x \hat y(\hat x +\hat q\hat y)(\hat x+\hat y/\hat q)]\,. 
\end{eqnarray*}
Hence $\la p\ra=\la xy(x+qy)(x+y/q) \ra$ for some $q\in [1/\sqrt 2,1)$. Next we show that polynomials of the form $xy(x+qy)(x+y/q)$ with different values of $q\in [1/\sqrt 2,1)$ belong to different orbits. For $A\in\GL(2,\RR)$, the projective transformation of $\RR P^1$ induced by $(A^{-1})^\top$  maps the set of zeroes of $p$ to the set of zeroes of $A \cdot p$. Each projective transformation of $\RR P^1$ preserves the cross-ratio of four points. Hence the set of all cross-ratios of the zeroes of a polynomial is an invariant of the $\GL(2,\RR)$-action. That is, if $z_1,\dots, z_4\in\RR P^1\cong \RR\cup\{\infty\}$ are the four different zeroes of $p\in P_4$, then
$$ C(\la p\ra) := \left\{ \left.\frac{z_i-z_k}{z_j-z_k}:\frac{z_i-z_l}{z_j-z_l}\ \right| \{i,j,k,l\}=\{1,2,3,4\}\right\}$$
is well defined. The zeroes of $xy(x+qy)(x+y/q)$ are $0,\infty,-q,-1/q$. Thus 
$$C(\la p\ra)=\{q^{\pm2}, (1-q^2)^{\pm1}, (1-1/q^2)^{\pm1} \}.$$
For $q\in(1/\sqrt 2,1)$, we have $C(\la p\ra)\cap(1/2,1)=q^2$. If $q=1/\sqrt 2$, then $C(\la p\ra)\cap(1/2,1)=\emptyset$. Thus different values of $q\in [1/\sqrt 2,1)$ give different orbits. Since $xy(x+qy)(x+y/q)=xy(x^2+rxy+y^2)$, where $r=q+1/q$,  we see that $\la p\ra=\la xy(x^2+rxy+y^2)\ra$ for exactly one $r\in (2,\,3/\sqrt 2]$.

(2) If $p$ has no zero, then $\la p\ra= \la p_1 p_2\ra$, where $p_1(x,y)=x^2+y^2$ and $p_2(x,y)=x^2+2bxy+cy^2$, $b^2<c$. The positive definite quadratic forms $p_1$ and $p_2$ are simultaneous diagonalisable, thus  $\la p\ra=\la (x^2+y^2)(x^2+sy^2)\ra$. Obviously, we can choose $s\in(0,1]$. 

(3) Suppose that $p$ has a  zero of multiplicity two and that all other zeroes (if further ones exist) are simple. 

(3)(a) If $p$ has two further zeroes, then  $\la p\ra=\la xy(ax+by)^2\ra$, $a\not=0,b\not=0$. Rescaling $x$ and $y$ we get 
$\la p\ra=\la xy(x+y)^2\ra=\la xy(x^2+rxy +y^2)\ra$ for $r=2$.

(3)(b) If $p$ has no further zero, then  $\la p\ra=\la x^2(x^2+bxy+cy^2)\ra=\la x^2(x^2+y^2)\ra$ by completing the square (with respect to $y$) and rescaling $x$ afterwards. In particular, $\la p\ra=\la (x^2+y^2)(x^2+sy^2)\ra$ for $s=0$.
\qed

Let $b_1,\dots,b_7$ be a basis as chosen in Lemma \ref{LII}. If $R\in {\mathcal K}(\fh)$, then 
$$R_{ij}:=R(b_i,b_j)=h(A^{ij},z^{ij},c^{ij}).$$
\begin{pr}\label{PholII} The space ${\mathcal K}(\fh)$ can be parametrised by real numbers $x_1,\dots,x_5, y_1,\dots,y_5, r_1,\dots,r_4,t,t_1,\dots,t_6, s_1, s_2, j_1,j_2$, where $R=h(A,z,c) $ $\in {\mathcal K}(\fh)$ is given by the data in Table 2. 
\begin{small}
\begin{table}
\renewcommand{\arraystretch}{1.5}
\begin{tabular}{|c|c|c|c|}
\hline
$R(b_i,b_j)$&$A$&$z$&$c$\\
\hline
$R_{16}$&$0$&$(x_4,x_3,x_2,x_1)$&$t_1+t$\\[0.5ex]
$R_{17}=\frac1{\sqrt2}R_{34}$&$0$&$(x_5,x_4,x_3,x_2)$&$t_4-t_5$\\[0.5ex]
$R_{26}=-\frac1{\sqrt2}R_{45}$&$0$&$(y_4,y_3,y_2,y_1)$&$t_2-t_3$\\[0.5ex]
$R_{27}$&$0$&$(y_5,y_4,y_3,y_2)$&$t_6+t$\\[0.5ex]
$R_{56}$&$\left(\begin{array}{cc} x_1&y_1 \\x_2& y_2\end{array}\right)$&$(t_6,t_2,s_2,j_2)$&$r_1$\\[3.5ex]
$R_{57}=\frac1{\sqrt2}R_{46}$&$\left(\begin{array}{cc} x_2&y_2 \\x_3& y_3\end{array}\right)$&$(t_5,t_1,t_3,s_2)$&$r_2$\\[3.5ex]
$R_{36}=\frac1{\sqrt2}R_{47}$&$\left(\begin{array}{cc} x_3&y_3 \\x_4& y_4\end{array}\right)$&$(s_1,t_4,t_1,t_2)$&$r_3$\\[3.5ex]
$R_{37}$&$\left(\begin{array}{cc} x_4&y_4 \\x_5& y_5\end{array}\right)$&$(j_1,s_1,t_5,t_6)$&$r_4$\\[3.5ex]
$R_{67}$&$\left(\begin{array}{cc} t_1+t&t_2-t_3 \\t_4-t_5& t_6+t\end{array}\right)$&$(r_4,r_3,r_2,r_1)$&$0$\\[3.5ex]
\hline
\multicolumn{4}{|c|}{$R_{12}=R_{13}=R_{14}=R_{15}=R_{23}=R_{24}=R_{25}=0$}\\[0.5ex]
\multicolumn{4}{|c|}{$R_{35}=R_{16}-R_{27}$}\\[0.5ex]
\hline
\end{tabular}
\\[0.5ex] 
\begin{center} 
    {{\bf  Table 2.}} 
\end{center}
\end{table}
\end{small}
\end{pr}
\proof 
Let $R$ be in ${\mathcal K} (\fh)$. Since $\la R_{ij}(b_k),b_l\ra=\la R_{kl}(b_i),b_j\ra$ and $R_{kl}\in\fh$, we have 
\begin{equation}\label{EhII1} R_{ij}=0,\quad i\le2\le j\le 5.
\end{equation}
The same argument gives 
$$
{\sqrt2} R_{17}=R_{34},\ {\sqrt2} R_{26}=-R_{45},\ {\sqrt2}R_{57}=R_{46},\ {\sqrt2} R_{36}=R_{47}
$$
and  
\begin{equation}\label{EhII2}
R_{35}=R_{16}-R_{27}.
\end{equation}
As above, we use the notation $b(i,j,k):=R_{ij}(b_k)+R_{jk}(b_i)+R_{ki}(b_j)$. 
From $b(i,j,6)=b(i,j,7)=0$, we conclude $A^{ij}=0$ for $i<j\le5$. 
By (\ref{EhII1}) and $b(i,4,6)=0$, we get $a_2^{i6}=0$ and $0=a_3^{i6}=a_1^{i7}=a_4^{i7}$, where the last identity follows from $b(i,3,7)=0$. Together with (\ref{EhII1}) and (\ref{EhII2}) this implies
$$A^{16}=A^{17}=A^{26}=A^{27}=0.$$
Let us consider the equations $b(i,j,k)=0$ for $i,j,k\in\{1,2,3,5,6,7\}$. These equations give, in particular,
$$\begin{array}{lll}
z_4^{67} = c^{56} =: r_1,&&z_2^{37}=z_1^{36}=:s_1,\\[0.5ex]
c^{57} = z_3^{67} =: r_2,&&z_3^{56}=z_4^{57}=:s_2,\\[0.5ex]
c^{36} = z_2^{67} =: r_3,&&z_1^{17}=a_3^{37}=:x_5,\\[0.5ex]
z_1^{67} = c^{37} =: r_4,&&z_4^{26}=a_2^{56}=:y_1.
\end{array}$$
The system of the remaining linear equations for the coefficients $A^{ij},\ z^{ij},\ c^{ij}$ that follow from $b(i,j,k)=0$ for $i,j,k\in\{1,2,3,5,6,7\}$ decomposes into six subsystems each of which is a system of equations in the elements of one of the following sets:
\begin{align*}
M_1:= &\{ a_2^{67},z_4^{35}, z_2^{56}, z_3^{57}, c^{26}\},\\ 
M_2:= &\{a_3^{67},z_2^{36}, z_3^{37}, z_1^{57}, c^{17}\},\\
M_3:= &\{a_1^{56},a_4^{56}, a_2^{57}, z_4^{16}, z_3^{26}, z_4^{27}, z_4^{35}, c^{25}\},\\
M_4:= &\{a_1^{67},a_4^{67}, z_3^{36}, z_4^{37}, z_1^{56}, z_2^{57}, c^{16}, c^{25}, c^{35}\},\\
M_5:= &\{a_2^{36},a_3^{56}, a_1^{57}, a_4^{57}, z_3^{16}, z_4^{17}, z_2^{26}, z_3^{27}, z_3^{35},c^{15}\},\\
M_6:= &\{a_1^{36},a_3^{36},a_4^{36},a_1^{37}, a_2^{37}, a_4^{37}, a_3^{57},\\ 
&\ \, z_1^{16},  z_2^{16},z_2^{17}, z_3^{17}, z_1^{26},  z_1^{27},z_2^{27}, z_1^{35}, z_2^{35}, c^{13}, c^{23}\}.
\end{align*}
The subsystem for $M_1$ is
$$z_4^{36}=z_2^{56}=:t_2,\quad z_2^{56}-a_2^{67}=z_3^{57}=:t_3,\quad c^{26}=a_2^{67}.$$
Similarly, $M_2$ is parametrised by $t_4$ and $t_5$ as claimed in the proposition. 
For $M_3$, we have
\begin{align*}
&a_1^{56}=z_4^{16}=:x_1, \quad a_2^{57}=a_4^{56}=:y_2, \quad z_3^{26}=z_4^{27},\\
&a_2^{57}-z_4^{27}= c^{25}=-a_4^{56}+z_3^{26},\quad z_4^{35}=a_1^{56}-a_4^{56},
\end{align*}
and for $M_4$,
\begin{align*}
&c^{27}=a_4^{67}=:t+t_6,\quad c^{16}=a_1^{67}=:t_1+t\\
&z_3^{36}-z_4^{37}=a_1^{67}-a_4^{67}=-z_1^{56}+z_2^{57},\\
&z_3^{36}-z_1^{56}=c^{35}=z_2^{57}-z_4^{37},
\end{align*}
where we put $t_1:=z_3^{36}$. The equations containing elements of $M_5$ are
\begin{align*}
&a_3^{56}=a_1^{57}=:x_2,\quad a_2^{36}=z_2^{26}=z_3^{27}=a_4^{57}=:y_3,\quad z_3^{16}=z_4^{17}, \\
&a_1^{57}-z_4^{17}=c^{15}=-a_3^{56}+z_3^{16},\\
&a_1^{57}-a_4^{57}=z_3^{35}=-a_2^{36}+a_3^{56}.
\end{align*}
Finally, for $M_6$, we get
\begin{align*}
&a_1^{36}=z_2^{16}=z_3^{17}=a_3^{57}=:x_3,\quad a_3^{36}=a_1^{37}=:x_4,\\
&z_1^{26}=z_2^{27}=:y_4,\quad a_4^{37}=z_1^{27}=:y_5\\
&z_1^{16}=z_2^{17},\quad a_4^{36}=a_2^{37},\\
&a_1^{37}-z_2^{17}=c^{13}=-a_3^{36}+z_1^{16},\quad a_2^{37}-z_2^{27}=c^{23}=-a_4^{36}+z_1^{26},\\
&a_1^{36}-a_4^{36}=z_2^{35}=-a_2^{37}+a_3^{57},\quad z_1^{35}=a_1^{37}-a_4^{37}.
\end{align*}

This shows that $R\in {\mathcal K}(\fh)$ is given as in Table 2. Moreover, the above considered systems of linear equations imply that each $R$ that is defined as in Table 2 for an arbitrary choice of the parameters $x_1,x_2,\dots, j_1,j_2$ satisfies $b(i,j,k)=0$ for $i,j,k\in\{1,2,3,5,6,7\}$. 
It is easy to show by a direct calculation that it also satisfies $b(i,j,4)=0$ for $1\le i,j\le7$.
\qed

The embedding of $\fgl(2,\RR)$ into $\fg_2^*$ defined by (\ref{EglII}) gives us an embedding of $\GL^+(2,\RR)$ into $\G$. If we send, moreover,  $\diag(1,-1)\in \GL(2,\RR)$ to $\diag(1,-1,-1,1,-1,1,-1)\in \G$, we obtain an embedding of $\GL(2,\RR)$ into $\G$, which we want to consider in this section. Note that this embedding is different from that defined by (\ref{EGL}).  With this identification we have
\begin{equation}\label{EAdg} 
\Ad(g)(h(A,z,c))= 
h(g A g^{-1},\, g\cdot z,\,\det(g) c).
\end{equation}
\begin{lm}
Either $\fa\in\{0,\ \fsl(2,\RR),\ \fgl(2,\RR)\}$ or the basis $b_1,\dots, b_7$ in Lemma \ref{LII} can be chosen such that $\fa$ is equal to one of the following Lie algebras:
\begin{enumerate}
\item $\RR\cdot C_a$, $\RR\cdot\diag(1,\mu),   \ \mu \in [-1,1]$;
\item $\fd$, $\fu(1)$.
\end{enumerate}
\end{lm}
\proof In a similar way as for Berger algebras of Type I we may conjugate $\fa$ by elements of $GL(2,\RR)$, now according to (\ref{EAdg}). Hence $\fa$ is one of the Lie algebras listed in Lemma~\ref{LaI}. Since $\fh$ acts semisimply on $S$, $\fa\cong\fgl(2,\RR)$ acts semisimply on $\RR^2$. This gives the assertion of the Lemma.  \qed

{\it Proof of Theorem \ref{T2}.} Below, we will use  the conjugation 
\begin{equation}\label{EAdII}
\Ad(\exp h(0,\bar z, 0))(h(A,z,c))=h(A,z-A\cdot\bar z,c-\eta(z,\bar z) -\textstyle{\frac12}\eta(\bar z, A\cdot \bar z))
\end{equation}
several times. We define $$Z:=\{z\in\RR^4\mid \exists \,c\in\RR:\ h(0,z,c)\in\fh\}.$$ Since $\fh$ is of Type II, we have $Z\not=0$. Obviously, $Z$ is invariant under $\fa$. Let $e_1,\dots,e_4$ be the standard basis of $\RR^4$ and denote by $Z(j_1,\dots,j_k)\subset \RR^4$, $k=1,2,3$, the span of $e_{j_1},\dots,e_{j_k}$.

(1) If $\fa\in \{\fsl(2,\RR),\, \fgl(2,\RR)\}$, then $Z=\RR^4$, since the action of $\fsl(2,\RR)$  on $P_d$  is irreducible. 

(2) Suppose $\fa\in\{\RR\cdot C_a,\fu(1)\}$. If $\fa=\RR\cdot C_a$, then $h(C_a,0,c_0)\in\fh$ for some $c_0\in\RR$ after a suitable conjugation of $\fh$ according to (\ref{EAdII}) since $C_a$ is bijective. If $\fa=\fu(1)$, then, possibly after conjugation, for all $U\in\fu(1)$, there is a real number $c$ such that $h(U,0,c)$ is in $\fh$. Indeed, $h(I,0,c_I)\in\fh$ after a suitable conjugation according to (\ref{EAdII}) since $I$ is bijective. Since 
$$[h(I,0,c_I),h(U,z_U,c_U)]=h(0,z_U,2c_U),$$
we have $h(U,0,-c_U)\in\fh$.
The restrictions of the representation of $\fgl(2,\RR)$ on $\RR^4$ to $\fu(1)$ and $\RR\cdot C_a$ decomposes into the two irreducible representations
$$Z_1:= \{(r,s,-r,-s)\mid r,s\in\RR\},\quad Z_2:= \{(3r,s,r,3s)\mid r,s\in\RR\}.$$
Indeed, both subspaces are invariant under $\fu(1)$ and they are irreducible since $C_a$ has eigenvalues $a\pm 3i$ on $Z_1$ and $a\pm i$ on $Z_2$. If $Z$ were equal to $Z_1$, then the non-isotropic vector $b_3+b_5$ would be in the kernel of $\fh$. Hence $Z=Z_2$ or $Z=\RR^4$, which gives the assertion.

(3) If $\fa=\fd$, then we can again conjugate $\fh$ according to (\ref{EAdII}) such that for all $D\in\fd$ there exists a $c\in\RR$ such that $h(D,0,c)\in\fh$. Indeed, as above $h(I,0,c_I)\in\fh$ after a suitable conjugation. Thus $[h(I,0,c_I),h(D,z_D,c_D)]=h(0,z_D,2c_D)\in\fh$ and $h(D,0,-c_D)\in\fh$ follows. The subspace $Z\subset \RR^4$ is invariant under $\fd$ if and only if it is invariant under $H$. Thus $Z$ is a direct sum of eigenspaces of $H$. Since $\fd$ considered as a subspace of $\fgl(2,\RR)$ is invariant under conjugation by $U:=\left(\begin{array}{cc} 0&-1\\1&0\end{array}\right)$, we may conjugate $\fh$ by $U\in\GL(2,\RR)\subset \G$. $U$ acts on $Z$ by $(z_1,z_2,z_3,z_4)\mapsto (z_4,-z_3,z_2,-z_1)$. Thus we may assume that $Z$ is one of the subspaces $\fn, Z(1), Z(2), Z(1,2), Z(1,3),$  $ Z(1,4),Z(2,3), Z(1,2,3),Z(1,2,4)$. For $Z=Z(1)$ and $Z=Z(1,4)$, $\fh$ is decomposable since $b_4$ is in the kernel. Moreover, we can exclude $Z=Z(2)$ and $Z=Z(1,2)$ since for these $Z$ the Lie algebra $\fh$ would be of Type III. For $Z(2,3), Z(1,2,3)$ and $Z(1,2,4)$ we get immediately $\fn(2,3)\subset\fh, \fn(1,2,3)\subset\fh$ and $\fn(1,2,4)\subset\fh$, respectively. Hence these three cases give Lie algebras that are on the list. For $Z=Z(1,3)$ we have to use Berger's criterion to show that $h(0,0,1)$ is in $\fh$.  The parameters of ${\mathcal K}(\fh)$ satisfy $x_1=\dots=x_5=y_1=\dots=y_4=t_1=\dots=t_6=s_1=s_2=0$. Hence $r_2\not=0$, which implies $h(0,0,1)\in\fh$. Thus $\fn(1,3)\subset\fh$, which gives the remaining Lie algebra in case (3). 

(4) Let $\fa$ be spanned by $\diag(1,\mu)$. The action of $\diag(1,\mu)$ on $\RR^4$ equals the multiplication by the matrix 
$D_0:=\diag(-1+2\mu, \mu, 1,2-\mu)$. Since $\mu\in[-1,1)$, all eigenvalues of $D_0$ are different. Hence every invariant subspace is spanned by elements of the standard basis of $\RR^4$.

(4)(a) Suppose first that $\mu\not\in\{0,1/2\}$. Then $D_0$ is invertible. Hence we may assume that $h_0:=h(\diag(1,\mu),0,c_0)\in\fh$ for some $c_0\in\RR$. Consequently, $\dim Z>1$ since otherwise $\fh$ would be decomposable. Furthermore, by indecomposability, $Z$ cannot be one of the spaces $Z(1,2)$ or $Z(1,4)$. Moreover,  $Z\not=Z(3,4)$ since $\fh$ is of Type II. The remaining spaces $Z(1,3)$ and $Z(2,4)$ can be excluded by Berger's criterion. Note first that the parameters $x_i$ and $y_i$ ($i=1,\dots,5$) of ${\mathcal K}(\fh)$ vanish since $\fa$ is spanned by $\diag(1,\mu)$ with $\mu\not=0$. The assumption $Z=Z(1,3)$ of $Z=Z(2,4)$ would imply $t_1=t_6=0$,  Thus $h_0\not\in \underline \fh$, which is a contradiction. 

(4)(b) Now suppose $\mu=1/2$. Then $(1,0,0,0)$ spans the kernel of $D_0$. Hence we may assume that $h(\diag(1,\mu),(t,0,0,0),c_0)\in\fh$ for certain $t,c_0\in\RR$. If $Z$ contains $Z(1)$, we choose $t=0$. If $t=0$, then we proceed as in (4)(a). Take now $t\not=0$. Then we can achieve $t=1$ conjugating by a suitable multiple of $I\in\GL(2,\RR)$. 
Again, $Z=Z(2)$ and $Z=Z(4)$ cannot occur because of in\-de\-com\-posability. For the remaining possibilities for $Z$ we want to use Berger's criterion.  In the same way as in (4)(a) we see that $x_i$ and $y_i$ vanish for $i=1,\dots,5$. Assume that $Z\subset Z(3,4)$ or $Z=Z(2,4)$. Then $t_1=t_6=0$, which as above leads to a contradiction. 

(4)(c) Finally, take $\mu=0$. The kernel of $D_0$ is spanned by $(0,1,0,0)$. Hence we may assume that $h(\diag(1,0),(0,t,0,0),c_0)\in\fh$ for certain $t,c_0\in\RR$. If $Z$ contains $Z(1)$, we choose $t=0$. Suppose first that $t=0$. As in (4)(a), $\dim Z>1$ and $Z\not\in\{Z(1,2), Z(1,4), Z(3,4)\}$. Moreover, $Z=Z(1,3)$ can be excluded by Berger's criterion. However, in contrast to the case $\mu\not=0$ we cannot rule out $Z=Z(2,4)$.  Indeed, if $Z=Z(2,4)$, then $\fh$ satisfies Berger's criterion if and only if $h(0,0,1)\in\fh$. Now we consider the case $t\not=0$. We may assume $t=1$. If $Z$ were in $\{Z(1), Z(2)\}$, then $\fh$ would be decomposable. For $Z\subset Z(1,3)$, $\fh$ would not satisfy Berger's criterion, see Proposition~\ref{PholII}.

(5) Now let $\fa$ be either trivial or equal to $\RR\cdot I$.  Then $Z$ can be an arbitrary subspace of $\RR^4$. We may conjugate $\fh$ by $\GL(2,\RR)\subset\G$ without changing  $\fa$. We want to use that in order to find a certain normal form for $Z$. Let $\la Z\ra$ denote the $\GL(2,\RR)$-orbit of $Z$. 

Let us first consider the case $\dim Z=1$. By Lemma~\ref{Lpol}, $Z$ is in one of the orbits $\la [\Phi_1^{-1}(x^3)]\ra$, $\la[ \Phi_1^{-1}(x^2y)]\ra$ or $\la[ \Phi_1^{-1}(x(x^2\pm y^2))]\ra$, where for $[z]$ denotes the line spanned by $z\in\RR^4$. For $[\Phi_1^{-1}(x^3)]=Z(4)$, the corresponding Lie algebra $\fh$ is decomposable since $b_4$ is in the kernel of $\fh$. For  $[\Phi_1^{-1}(x^2y)]=Z(3)$, we get a Lie algebra $\fh$ of Type III. For $[\Phi_1^{-1}(x(x^2\pm y^2))]=\{(0,z/3,0,z)\mid z\in\RR\}$, the corresponding Lie algebra $\fh$ is decomposable since $3b_3\mp b_5$ is in the kernel.

Now we turn to the case $\dim Z=2$. Let us consider the Pl\"ucker embedding of the Grassmannian of 2-planes in $\RR^4$ into $\PP(\bigwedge^2\RR^4)$. The line spanned by 
$$t_0w_0+z_{12} e_{12} +z_{13} e_{13}+t'w'+z_{24} e_{24}+z_{34} e_{34}\in\textstyle{\bigwedge^2\RR^4}$$ 
is in the image of this embedding if and only if 
\begin{equation}\label{EPl}
z_{12}z_{34}-z_{13}z_{24}= 3(t_0^2-t'^2)\,.
\end{equation}
Let us denote by $\hat Z$ the image of $Z$ under the Pl\"ucker embedding and by $\la\hat Z\ra$ the $\GL(2,\RR)$-orbit of $\hat Z$. Furthermore, for $\alpha\in \bigwedge^2\RR^4$ denote by $[\alpha]$ the element of $\PP(\bigwedge^2\RR^4)$ represented by  $\alpha$. Then $\la \hat Z\ra=\la [\Phi_2^{-1}(p)+t_0w_0]\ra$ for some polynomial $p\in P_4$ and some $t_0\in\RR$. By Lemma~\ref{Lpol} we may assume that $p$ is one of the polynomials $x^4,\, x^3y,\, x^2y^2,\, xy(x^2+rxy+y^2)$ for $r\in[0,\,3/\sqrt 2]$ or $(x^2+y^2)(x^2+sy^2)$ for $s\in[0,1]$. Let us start with $p=x^4$. Since $\Phi_2^{-1}(x^4)+t_0w_0=e_{34}+t_0w_0$ is in the image of the Pl\"ucker embedding if and only if $t_0=0$, we see that $\hat Z$ is in the $\GL(2,\RR)$-orbit of $e_{34}$. Since the Pl\"ucker embedding is $\GL(2,\RR)$-equivariant, $Z$ is in the $\GL(2,\RR)$-orbit of $Z(3,4)$. But then $\fh$ is of Type III. Now we consider $p=x^3y$. We have $\la \hat Z\ra= \la[\Phi_2^{-1}(x^3y)+t_0w_0]\ra= \la[\Phi_2^{-1}(y^3x)+t'_0w_0]\ra=\la [e_{13}+t''_0w_0]\ra$. Since $e_{13}+t''_0w_0$ has to be in the image of the Pl\"ucker imbedding, we get $t''_0=0$. Thus $Z$ is in the $\GL(2,\RR)$-orbit of $Z(1,3)$. In the same way as in case (3) we can show that $h(0,0,1)$ is in $\fh$ using Berger's Criterion. Hence $\fh=\fa\ltimes \fn(1,3)$ up to conjugation. For $p=x^2y^2$, we get $\la \hat Z\ra= \la[\Phi_2^{-1}(x^2y^2)+t_0w_0]\ra= \la [w'+t'_0w_0]\ra$. By (\ref{EPl}), $t'_0=\pm1$. If $t_0'=1$, then $\la \hat Z\ra= \la w'+w_0\ra=\la e_{23}\ra$. Hence $Z$ is in the orbit of $Z(2,3)$, which implies $\fh=\fa\ltimes \fn(2,3)$. For $t'_0=-1$ we obtain
$\la \hat Z\ra= \la w'-w_0\ra=\la e_{14}\ra$. Thus $Z$ is in the orbit of $Z(1,4)$. But then $\fh$ is decomposable since $b_4$ is in the kernel of $\fh$. Now take $p=xy(x^2+rxy+y^2)$, $r\in[0,\,3/\sqrt 2]$.  Then 
\begin{align}\label{EZ1}
\la \hat Z\ra &= \textstyle {\la[\Phi_2^{-1}(p)+t'_0w_0]\ra= \la [\frac12 e_{24}+\frac r6 w' +\frac12 e_{13} +t'_0w_0]\ra}\nonumber\\
&=\textstyle{\la [e_{24} +\frac r3 w' + e_{13} +t_0w_0]\ra}\nonumber\\
&=\textstyle {\la[e_{13}+(\frac r3+t_0)e_{23}+(r-3t_0)e_{14}+e_{24}]\ra.} 
\end{align}
By (\ref{EPl}), we have $9t_0^2=r^2-3$, which is possible only for $r\ge\sqrt 3$. The map
\begin{eqnarray*}
&M:= \left\{(r,t_0)\left|\ 9t_0^2=r^2-3,\ r\in[\sqrt3,\,3/\sqrt 2]\right.\right\} \longrightarrow  \textstyle {\left[\frac{\sqrt 3-1}{\sqrt6},\frac{\sqrt3+1}{\sqrt 6}\right]}&\\
&(r,t_0)\longmapsto \alpha(r,t_0):=\textstyle{\frac r3}+t_0&
\end{eqnarray*}
is a bijection. Note that $\alpha(r,t_0)\cdot(r-3t_0)=1$. Now (\ref{EZ1}) implies
 \begin{align*}
 \la \hat Z\ra&=\la[e_{13}+\alpha e_{23}+(1/\alpha) e_{14}+e_{24}]\ra=\la [\alpha e_{13}+\alpha^2e_{23}+e_{14}+\alpha e_{24}]\ra\\
 &=\la (e_1+\alpha e_2)\wedge(\alpha e_3 +e_4)]\ra\end{align*}
for $\alpha:=\alpha(r,t_0)$.
Hence $Z$ and $\{(z_1,\alpha z_1, \alpha z_4, z_4)\mid z_1, z_4\in\RR\}$ are in the same orbit. Thus we are in case (\ref{Ic}) of the theorem. At last, we consider $p=(x^2+y^2)(x^2+sy^2)$, $s\in[0,1]$.  Then
\begin{align}\label{EZ2}
\la \hat Z\ra &= \textstyle {\la[\Phi_2^{-1}(p)+t_0w_0]\ra= \la [e_{34}+\frac {s+1}6 w' +s e_{13} +t_0w_0]\ra}\nonumber\\
&=\textstyle{\la [e_{24} +\frac r3 w' + e_{13} +t_0w_0]\ra}\nonumber\\
&=\textstyle {\la[s e_{12}+3(\frac {s+1}6-t_0)e_{14}+(\frac {s+1}6 +t_0)e_{23}+e_{34}]\ra.} 
\end{align}
By (\ref{EPl}), we have $s=3\left(t_0^2-\left(\frac
{s+1}6\right)^2\right)$. Thus, for a given parameter $s$, the coefficients $\frac {s+1}6-t_0$ and $\frac {s+1}6+t_0$ appearing in (\ref{EZ2}) are exactly the roots $\alpha_{1,2}$ of $3\alpha^2-(s+1)\alpha-s$. In particular, $\alpha_1\alpha_2=-s/3$. Now suppose in addition that $s\not=0$. Then the latter equation implies $3\alpha_1=-s/\alpha_2$. Consequently, 
\begin{align*}
\la \hat Z\ra &= \textstyle  {\la[s e_{12}-(s/\alpha)e_{14}+\alpha e_{23}+e_{34}]\ra} 
= \textstyle  {\la[s \alpha e_{12}-s e_{14}+\alpha^2 e_{23}+\alpha e_{34}]\ra} \\
&= \la [\, (se_1-\alpha e_3)\wedge (\alpha e_2-e_4)\,]\ra
\end{align*}
where $\alpha\in\{\alpha_1,\alpha_2\}$. Hence, if $s\not=0$, then $Z$ is in the $\GL(2,\RR)$-orbit of $\{(sz_1,\alpha z_2, -\alpha z_1, -z_2)\mid z_1, z_2\in \RR\}$, which implies that $\fh$ is conjugated to the Lie algebra in (\ref{Id}). If $s=0$, then (\ref{EZ2}) implies that $\la \hat Z\ra=\la[e_{14}+e_{34}]\ra=\la[(e_1+e_3)\wedge e_4]\ra$ or
$$\textstyle{\la \hat Z\ra=\la[\frac 16 (w'+w_0)+e_{34}]\ra=\la[-\frac 12 (w'+w_0)- e_{34}]\ra=\la[e_3\wedge(e_2-e_4)]\ra}.$$ 
Thus we are in case (\ref{Ia}) or (\ref{Ib}), respectively.

Finally, we study the case $\dim Z=3$. Recall that the action of $\GL(2,\RR)$ on  $\bigwedge^3\RR^4$  is equivalent to the representation on $P_3$, see~(\ref{Eaequ}). By Lemma~\ref{Lpol}, the image of $Z$ under the Pl\"ucker embedding is in the same $\GL(2,\RR)$-orbit as either $[\Phi_3^{-1}(x^3)]$, $[\Phi_3^{-1}(x^2y)]$ or $[\Phi_3^{-1}(x(x^2\pm y^2))]$. Because of $\Phi_3^{-1}(x^3)=e_{234}$, $\Phi_3^{-1}(x^2y)=e_{134}$ and $[\Phi_3^{-1}(x(x^2\pm y^2))]=[e_{234}\pm e_{124}]=[(e_1\pm e_3)\wedge e_2\wedge e_4]$ we see that $\fh$ is conjugated to $\fa\ltimes \fn_1$, where $\fn_1=\fn(2,3,4)$ or $\fn_1=\fn(1,3,4)$ or $\fh$ is conjugated to the Lie algebra in item (\ref{Ie}).
\qed    

\section{Holonomy of symmetric spaces admitting a $\G$-structure} \label{S3}

Indefinite symmetric spaces of signature $(4,3)$ with a $\G$-structure were classified in \cite{Ka}. Their holonomy algebras can be easily read off from this classification. They are abelian and two- or three-dimensional (unfortunately, there is an obvious mistake in the formulation of Cor.~6.9 in \cite{Ka}). Let us check how they fit into the classification of holonomy algebras in Section~\ref{S2}. 

Let $X=G/G_+$ be a (pseudo-)Riemannian symmetric space, where $G$ is the transvection group of $X$. The reflection at the base point $x_0:=eG_+\in G/G_+$ defines an involution $\theta$ on the Lie algebra $\fg$ of $G$. Let $\fg_+$ and $\fg_-$ denote the eigenspaces of $\theta$ with eigenvalue $1$ and $-1$, respectively. Then $\fg_-$ can be identified with the tangent space of $X$ at $x_0$ and $\fg_+$ can be identified with the holonomy algebra (as an abstract Lie algebra). The holonomy representation is given by the adjoint representation of $\fg_+$ on $\fg_-$.  

The classification of symmetric spaces with $\G$-structure is given by the list in Theorem 6.8 in \cite{Ka}. Item 1 of this list contains a one-parameter family of symmetric spaces and each of the items 2 (a) and (b) contains a single space. 

Let us first consider the family in item 1. In the notation of the theorem, the holonomy algebra $\fg_+$ is spanned by $Z_B,A_1,B$ and the tangent space $\fg_-$ by $Z_1,Z_2,Z_3,A,L_1,L_2,L_3$. The adjoint representation of $\fg_+$ on $\fg_-$ is defined by the Lie bracket given by Equation~(8) in~\cite{Ka}. It is easy to verify that the holonomy algebra as a subalgebra of $\fso(4,3)$ is of Type II and equals $\fn(1,3)$ with respect to the basis $Z_2,Z_3,L_1,A,Z_1,L_2,L_3$. 

The space in item 2\,(a) has a holonomy algebra $\fg_+=$ spanned by $Z_B,A_1,B$. As a subalgebra of $\fso(4,3)$ the holonomy algebra is of Type III and equals $\fm(1,0,2)$ with respect to the basis $Z_3,Z_2,Z_1,A,L_3,L_2,L_1$. For the symmetric space in item 1\,(b) we get a similar result. Its holonomy algebra is spanned by $Z_B,B$ and using the adjoint representation it can be identified with $\fm(1,0,1)$ with respect to the same basis of $\fg_-$. In particular, it is also of Type III.
\section{Left-invariant metrics with holonomy in $G^*_2$ }\label{S4}
In this section we construct examples of Lie groups $G$ admitting left-invariant torsion free $G_2^*$-structures whose holonomy algebra is one of the list in Theorems \ref{T1}, \ref{T3} and \ref{T2}. The $G_2^*$-structure on $G$ is induced by a three-form $\omega$ on the Lie algebra $\frak g$ of $G$. If we denote by $\ip$ the induced inner product, the Levi-Civita connection $\nabla$  is  determined by
$$
2\langle \nabla_u v, w \rangle =   \langle [x,y], z \rangle - \langle [y,z], x \rangle + \langle [z,x], y \rangle , \quad x,y,z \in \frak g.
$$
If $b_1, \ldots, b_7$ is a basis of $\frak g$ we will denote by $\Lambda_j$ the endomorphism corresponding to $\nabla_{b_j}$ and by $R_{kl}$ the curvature operator $R(b_k ,b_l)$. We compute the holonomy algebra  using the Ambrose-Singer holonomy theorem \cite{AS}, which states that the holonomy  algebra  is spanned by the curvature operators $R_{xy}, x,y  \in {\frak g},$ together with their covariant derivatives.

\subsection{Examples of Type I}

In all  examples of this subsection, $b_1,\dots,b_7$ will be a basis of a Lie algebra $\fg$ such that $\omega$ equals $$\sqrt2 (b^{127}+b^{356})-e^4\wedge(b^{15}+b^{26}-b^{37})$$ and the induced inner product on $\fg$ equals 
$$2(b^1\cdot e^5+b^2\cdot e^6+b^3\cdot e^7)- (b^4)^2.$$ 

\begin{ex}   (5-dimensional holonomy) Let $\frak g$ be  the  Lie algebra with structure equations
$$\begin{array}{l}
db^1 = -2b^{15} - b^{56}, \\[2pt]
db^2 = -2b^{25}-b^{35}-b^{56}, \\[2pt]
db^3 = -b^{35}, \\[2pt]
db^4 = b^{45}-\sqrt2 b^{67}, \\[2pt]
db^5 = 0,\\[2pt]
db^6 = - b^{56}, \\[2pt]
db^7 = b^{36}- b^{56}. 
\end{array}$$
Then
\begin{align*}
\Lambda_1 &= 0\\
\Lambda_2&= \textstyle{\frac12}(b^5_2-b^6_1)\\
\Lambda_3 &= b^3_2-b^6_7+\textstyle{\frac12}(b^5_3-b^7_1),\\
\Lambda_4 &=-b^4_1-b^5_4+\textstyle{\frac1{\sqrt2}}(b^6_3-b^7_2),\\
\Lambda_5 &=-2(b^1_1-b^5_5)-\textstyle{\frac32}(b^2_2-b^6_6)-\textstyle{\frac12}(b^3_3-b^7_7)-b^3_2+b^6_7-b^5_2+b^6_1,\\
\Lambda_6 &= -\textstyle{\frac12}(b^2_1-b^5_6)+\textstyle{\frac1{\sqrt2}}(b^4_3+b^7_4)-b^5_2+b^6_1,\\
\Lambda_7 &=  -\textstyle{\frac12}(b^3_1-b^5_7)-\textstyle{\frac1{\sqrt2}}(b^4_2+b^6_4).\end{align*}
For the holonomy algebra $\fh$, we obtain
$$\fh=\Span\{R_{25},R_{35}, R_{45},R_{56}, R_{57}\}=\fm.$$
\end{ex}
\begin{ex} (6-dimensional holonomy)  Let $\frak g$ be  the  Lie algebra with structure equations
$$
\begin{array}{l}
d b^1 = d b^5 = d b^6 =0,\\[2pt]
d b^2 = \sqrt2 (b^{25} + b^{35} - b^{56} - b^{57}),\\[2pt]
d b^3= -db^7 = -\frac1{\sqrt2}b^{35} - b^{46} - \sqrt2 b^{56} +\frac1{\sqrt2}b^{57},\\[2pt]
db^4 = - b^{36} + b^{67}.
\end{array}$$
Then
\begin{align*}
\Lambda_1 &= 0\\
\Lambda_2&= -\textstyle{\frac1{\sqrt2}}(b^5_2-b^6_1),\\
\Lambda_3 &= \textstyle{\frac1{\sqrt2}}(b^3_1-b^5_7)+b^4_2+b^6_4,\\
\Lambda_4 &=b^3_2-b^6_7,\\
\Lambda_5 &=\textstyle{\frac1{\sqrt2}}(b^2_2-b^6_6-b^3_3+b^7_7)+\sqrt2(b^3_2-b^6_7),\\
\Lambda_6 &= \textstyle{\frac1{\sqrt2}}(b^2_1-b^5_6)-b^4_3-b^7_4-{\sqrt2}(b^5_2-b^6_1+b^5_3-b^7_1),\\
\Lambda_7 &= - \sqrt2(b^5_2-b^6_1)+\textstyle{\frac1{\sqrt2}}(b^5_3-b^7_1).\end{align*}
For the holonomy algebra $\fh$, we obtain
$$\fh=\Span\{R_{25},R_{35},R_{36}, R_{56},R_{67},\nabla_{e_6}R_{36}\}=\RR\cdot N\ltimes\fm.$$
\end{ex}
\begin{ex}   (7-dimensional holonomy) Let $\frak g$ be  the  Lie algebra with structure equations
$$
\begin{array}{l}
d b^1 = -   b^{16} -  \frac 1{ \sqrt{2} }b^{26} + \frac12 b^{36} + (1+\sqrt2) b^{46} +b^{56} -b^{67},\\[2pt]
d b^2 = -\frac32 b^{26}+(\frac{\sqrt2}4+\frac12)b^{36}+(3+\frac1{\sqrt2})b^{46}+(\frac1{\sqrt2}-1)  b^{56}-(\frac1{\sqrt2}+1)b^{67},\\[2pt]
d b^3 =  \frac{\sqrt2}4 b^{35} + \frac 12 b^{36} + \frac12 b^{46} - b^{56} - \frac1{\sqrt2} b^{57},\\[2pt]
de^4 = \frac 14  b^{36} +  (1+ \sqrt{2}) b^{56} - \frac 12 b^{67}  ,\\[2pt]
d b^5= b^{56},\\[2pt]
d b^6 =0,\\[2pt]
d b^7 =  - \frac{ \sqrt{2}}8 b^{35}  - \frac 12 b^{36} - \frac 14 b^{46} -(\frac12+3\sqrt2)b^{56}+ \frac{ \sqrt{2}}4 b^{57} + \frac 12 b^{67}.
\end{array}
$$
Then
\begin{align*}
\Lambda_1 &= 0\\
\Lambda_2&= -\textstyle{\frac{\sqrt2}4}(b^5_2-b^6_1),\\
\Lambda_3 &= -\textstyle{\frac{\sqrt2}8}(b^3_1-b^5_7)-\textstyle{\frac12}(b^3_2-b^6_7)-\textstyle{\frac14}(b^4_2+b^6_4) -\textstyle{\frac3{\sqrt2}}(b^5_2-b^6_1),\\
\Lambda_4 &=-\frac14(b^3_2-b^6_7),\\
\Lambda_5 &=\textstyle{\frac{\sqrt2}4}(-b^2_2+b^6_6+b^3_3-b^7_7)-\textstyle{\frac3{\sqrt2}}(b^3_2-b^6_7)+b^5_2-b^6_1,\\
\Lambda_6 &= -b^1_1+b^5_5-\textstyle{\frac32}(b^2_2-b^6_6) +\textstyle{\frac12}(b^3_3-b^7_7) -\textstyle{\frac{\sqrt2}4}(b^2_1-b^5_6)  +(\frac12+\frac3{\sqrt2})(b^3_1-b^5_7)
\\&\quad +(  \textstyle{\frac{\sqrt2}4+\frac12})(b^3_2-b^6_7)+(1+\sqrt2)(b^4_1+b^5_4)+(3+\textstyle{\frac1{\sqrt2}})(b^4_2+b^6_4) \\
&\quad +\textstyle{\frac12}(b^4_3+b^7_4)+(\textstyle{\frac1{\sqrt2}}-1)(b^5_2-b^6_1)-b^5_3+b^7_1) -(\textstyle{\frac1{\sqrt2}}+1)(b^6_3-b^7_2),
\\
\Lambda_7 &=  -\textstyle{\frac1{\sqrt2}}(b^5_3-b^7_1).\end{align*}
For the holonomy algebra $\fh$, we obtain
$$\fh=\Span\{R_{25}, R_{35}, R_{36}, R_{45}, R_{56}, R_{57}, \nabla_{b_5}R_{56}\}=\fs_{1/2}\ltimes\fm.$$
\end{ex}
\subsection{Examples of Type II}
In all  examples of this subsection, $b_1,\dots,b_7$ will be a basis of a Lie algebra $\fg$ such that $\omega$ equals $$\sqrt 2(-b^{157}+b^{236})-b^4\wedge(b^{16}-b^{27}-b^{35})$$ and  the  induced inner product on $\fg$ equals 
$$2(b^1\cdot b^6+b^2\cdot b^7+b^3\cdot b^5)- (b^4)^2.$$ 
\begin{ex}   (3-dimensional holonomy) Let $\frak g$ be  the  Lie algebra with structure equations
$$
\begin{array}{l}
d b^1 =  - b^{17}+\frac1{\sqrt2}b^{34} +( 1 - \sqrt{2}) b^{37}  - b^{46} -  b^{57} + b^{67},\\[4pt]
d b^2 =  - b^{27}+b^{37}  - b^{47} - b^{67}\\[2pt]
d b^3 = (1 - \sqrt{2})  b^{67},\\[2pt]
db^4 =  - b^{47},\\[2pt]
d b^5= -b^{37} +\frac1{\sqrt{2}} b^{46} + b^{67},\\[2pt]
d b^6 = b^{67},\\[2pt]
d b^7=  0.
\end{array}
$$
Then 
\begin{align*}
\Lambda_1 &=\Lambda_2=0,\\
\Lambda_3 &= -b^3_2 + b^7_5 + (1-\textstyle{\frac1{\sqrt 2}} )(b^6_2-b^7_1),\\
\Lambda_4 &= \textstyle{\frac1{\sqrt 2}} (b^3_1-b^6_5)+(b^4_2+e^7_4),\\
\Lambda_5 &=\textstyle{\frac1{\sqrt 2}}(b^7_1-b^6_2),\\
\Lambda_6 &= (1-\textstyle{\frac1{\sqrt 2}})(b^3_2-b^7_5)-\textstyle{\frac1{\sqrt 2}}(b^5_2-b^7_3)+b^6_2-b^7_1-b^4_1-b^6_4,\\
\Lambda_7 &= b^3_2-b^7_5-\textstyle{\frac1{\sqrt 2}}(b^3_1-b^6_5)-b^2_2+b^7_7-b^1_1+b^6_6-b^4_2-b^7_4  \\ 
&\quad -(1-\textstyle{\frac1{\sqrt 2}})(b^5_1-b^6_3)-b^6_2 +b^7_1.
\end{align*}
For the holonomy algebra $\fh$, we obtain
$$
\fh=\Span\{ R_{37}, R_{46}, R_{67}\}=\fn(1,3).$$
\end{ex}
\begin{ex}   (5-dimensional holonomy) Let $\frak g$ be  the  Lie algebra with structure equations
$$\begin{array}{l}
db^1 = 2b^{15}+4b^{56}+b^{57}, \\[2pt]
db^2 =b^{17}+b^{25}+\sqrt2 b^{34}-\sqrt2 b^{46}-e^{56}-b^{57}-b^{67}, \\[2pt]
db^3 = b^{35}-3b^{56}, \\[2pt]
db^4 = -\sqrt2 b^{37}-\sqrt2 b^{67}, \\[2pt]
db^5 = 0,\\[2pt]
db^6 = 2 b^{56}, \\[2pt]
db^7 =  b^{57}. 
\end{array}$$
Then
\begin{align*}
\Lambda_1 &= \Lambda_2=\Lambda_3 =0,\\
\Lambda_4 &= \sqrt 2 (b^3_2-b^7_5+b^6_2-b^7_1),\\
\Lambda_5 &=2(b^1_1-b^6_6)+b^2_2-b^7_7+b^3_3-b^5_5-3(b^5_1-b^6_3)+b^6_2-b^7_1,\\
\Lambda_6 &= 4(b^5_1-b^6_3),\\
\Lambda_7 &= b^1_2-b^7_6-\sqrt 2(b^3_4+b^4_1+b^4_5+b^6_4)-b^5_2+b^7_3-b^6_2+b^7_1.\end{align*}
For the holonomy algebra $\fh$, we obtain
$$\fh=\Span\{R_{37},R_{46}, R_{56},R_{67}, \nabla_{b_7}R_{67}\}=\fn.$$
\end{ex}
\begin{ex}   (8-dimensional holonomy) Let $\frak g$ be  the  Lie algebra with structure equations
$$\begin{array}{l}
db^1 = -b^{17}+b^{26}-b^{36}-b^{37}-b^{46}+e^{56}-b^{57}-b^{67}, \\[2pt]
db^2 =\textstyle{\frac 53b^{16} -b^{27}+\frac{13}3b^{36}+\frac{11}9b^{37}-\frac43 \sqrt2 b^{46}+(\frac83 \sqrt2 +\frac53) b^{47} }\\ [2pt]
\quad\quad\quad \textstyle{-9b^{56}-\frac13 b^{57}-b^{67}} \\[2pt]
db^3 = b^{37}-2\sqrt2 b^{46}-b^{57}-(\sqrt2-10)b^{67}, \\[2pt]
db^4 = \textstyle{\frac23 b^{47}+\frac23\sqrt2 b^{67}, }\\[2pt]
db^5 = \textstyle{\frac79 b^{37}-\frac23 \sqrt2 b^{46}-\frac53 b^{57}-(2+\frac53 \sqrt2)b^{67},}\\[2pt]
db^6 = db^7 =  0. 
\end{array}$$
Then
\begin{align*}
\Lambda_1 &= \textstyle{-\frac43(b^6_2-b^7_1)}\\
\Lambda_2&=0\\
\Lambda_3 &=\textstyle{\frac79(b^3_2-b^7_5)-\frac{\sqrt2}3(b^4_1+b^6_4)-\frac13(b^5_2-b^7_3)-(\frac{11}3+\frac56\sqrt2)(b^6_2-b^7_1)},\\
\Lambda_4 &= \textstyle{-\frac{\sqrt2}3(b^3_1-b^6_5)-\frac23(b^4_2+b^7_4)-\sqrt2(b^5_1-b^6_3)+\frac{\sqrt2}3(b^6_2-b^7_1)},\\
\Lambda_5 &=\textstyle{-\frac13(b^3_2-b^7_5)-\sqrt2(b^4_1+b^6_4)-b^5_2+b^7_3+(9-\frac1{\sqrt2})(b^6_2-b^7_1)}\\
\Lambda_6 &= \textstyle{\frac13(b^1_2-b^7_6)+b^2_1-b^6_7 -b^3_1+b^6_5 + (\frac23 -\frac56\sqrt2)(b^3_2-b^7_5) -\frac{\sqrt2}3 (b^3_4+b^4_5) }\\
&\textstyle{\quad -b^4_1-b^6_4 -\sqrt2(b^4_2+b^7_4+b^4_3+b^5_4) +b^5_1-b^6_3 -\frac1{\sqrt2}(b^5_2-b^7_3)-b^6_2+b^7_1},\\
\Lambda_7 &=\textstyle{\frac13(b^1_1-b^6_6)-b^2_2+b^7_7+\frac43(b^3_3-b^5_5)+(\frac83+\frac56\sqrt2)(b^3_1-b^6_5)+\frac{11}9(b^3_2-b^7_5)}\\
&\quad \textstyle{-\frac13\sqrt2(b^4_1+b^6_4)+(\frac83\sqrt2+\frac53)(b^4_2+b^7_4)-(10-\frac1{\sqrt2})(b^5_1-b^6_3)}\\
&\quad\textstyle{-\frac13(b^5_2-b^7_3)-b^6_2+b^1_7   }.
\end{align*}
For the holonomy algebra $\fh$, we obtain
$$\fh=\Span\{R_{17}, R_{36}, R_{37}, R_{56}, R_{57},R_{67}, \nabla_{b_6}R_{67}, \nabla_{b_7}R_{67}\}=\fsl(2,\RR)\ltimes\fn.
$$
\end{ex}
\subsection{Example of Type III}
In the   example of this subsection, $b_1,\dots,b_7$ will be a basis of a Lie algebra $\fg$ such that $\omega$ equals $$\sqrt2 (b^{127}+b^{356})-e^4\wedge(b^{15}+b^{26}-b^{37})$$ and the induced inner product on $\fg$ equals 
$$2(b^1\cdot e^5+b^2\cdot e^6+b^3\cdot e^7)- (b^4)^2.$$ 

\begin{ex}   (3-dimensional holonomy) Let $\frak g$ be  the  Lie algebra with structure equations
$$\begin{array}{l}
db^1 = -b^{15}-b^{45}-b^{56}, \\[2pt]
db^2 = -\textstyle{\frac13}b^{25}-b^{35}-b^{36}-\sqrt2 b^{45}-b^{56}+(\sqrt2 -1)b^{57}+b^{67}, \\[2pt]
db^3 = -\textstyle{\frac23}b^{35}-b^{56}+\textstyle{\frac43}b^{57}, \\[2pt]
db^4 = -\sqrt2 b^{56}, \\[2pt]
db^5 = 0,\\[2pt]
db^6 = -\frac13 b^{56}, \\[2pt]
db^7 = -b^{56}-\frac23 b^{57}. 
\end{array}$$
Then
\begin{align*}
\Lambda_1 &= \Lambda_2=\Lambda_3 =0\\
\Lambda_4 &= \sqrt 2 (b^5_2-b^6_1)\\
\Lambda_5 &=-b^1_1+b^5_5-\textstyle{\frac13}(b^2_2-b^6_6)-\textstyle{\frac23}(b^3_3-b^7_7)-b^3_2+b^6_7-b^4_1-b^5_4-b^5_2+b^6_1\\&\quad +\textstyle{\frac1{\sqrt 2}}(b^6_3-b^7_2)\\
\Lambda_6 &= -b^3_2+b^6_7-\sqrt2 (b^4_1+b^5_4)-b^5_2+b^6_1 -(1-\textstyle{\frac1{\sqrt 2}})(b^5_3-b^7_1) +b^6_3-b^7_2\\
\Lambda_7 &= -(1-\textstyle{\frac1{\sqrt 2}})(b^5_2-b^6_1)+\textstyle{\frac43}(b^5_3-b^7_1).\end{align*}
For the holonomy algebra $\fh$, we obtain
$$\fh=\Span\{R_{45}, R_{56},R_{57}\}=\fm(1,0,2).$$
\end{ex}

\bigskip

\bigskip

\small\noindent Dipartimento di Matematica G. Peano, Universit\`a di
Torino, Via Carlo Alberto 10, Torino, Italy.\\
\texttt{annamaria.fino@unito.it}

\small\noindent Institut f\"ur Mathematik und Informatik, Universit\"at Greifswald,
Walther-Rathenau-Str. 47,  D-17487 Greifswald, Germany. \\
\texttt{ines.kath@uni-greifswald.de}

\end{document}